# The Quadrahelix:
# A Nearly Perfect Loop of Tetrahedra

Michael Elgersma (Minneapolis, MN, USA); elgersma.michael@gmail.com
Stan Wagon (Macalester College, St. Paul, MN, USA); wagon@macalester.edu

Dedicated to the memory of Stash Świerczkowski (1932–2015), mathematician and adventurer.

**Abstract.** In 1958, S. Świerczkowski proved that there cannot be a closed loop of congruent interior-disjoint regular tetrahedra that meet face-to-face. Such closed loops do exist for the other four regular polyhedra. It has been conjectured that, for any positive $\epsilon$, there is a tetrahedral loop such that its difference from a closed loop is less than $\epsilon$. We prove this conjecture by presenting a very simple pattern that can generate loops of tetrahedra in a rhomboid shape having arbitrarily small gap. Moreover, computations provide explicit examples where the error is under $10^{-100}$ or $10^{-10^6}$. The explicit examples arise from a certain Diophantine relation whose solutions can be found through continued fractions; for more complicated patterns a lattice reduction technique is needed.

## 1. Introduction

For each of the Platonic solids except the regular tetrahedron, it is easy to construct embedded face-to-face chains using congruent copies; *embedded* means that no two polyhedra have nonempty interior intersection. Figure 1 shows such toroidal loops of length 8 for octahedra, dodecahedra, and icosahedra; cubes are trivial. Steinhaus [13] wondered whether such loops exist for tetrahedra and in 1958 S. Świerczkowski [14] provided the answer: there is no such tetrahedral chain (embedded or not). We present the proof found by Dekker [5] in 1959 and, independently, Mason [11] in 1972. More details are in [16]. The main point is that the group generated by the reflections in the four tetrahedral faces is the free product $\mathbb{Z}_2 * \mathbb{Z}_2 * \mathbb{Z}_2 * \mathbb{Z}_2$. We use $I$ for the identity matrix; if the dimension is not clear a subscript is used.

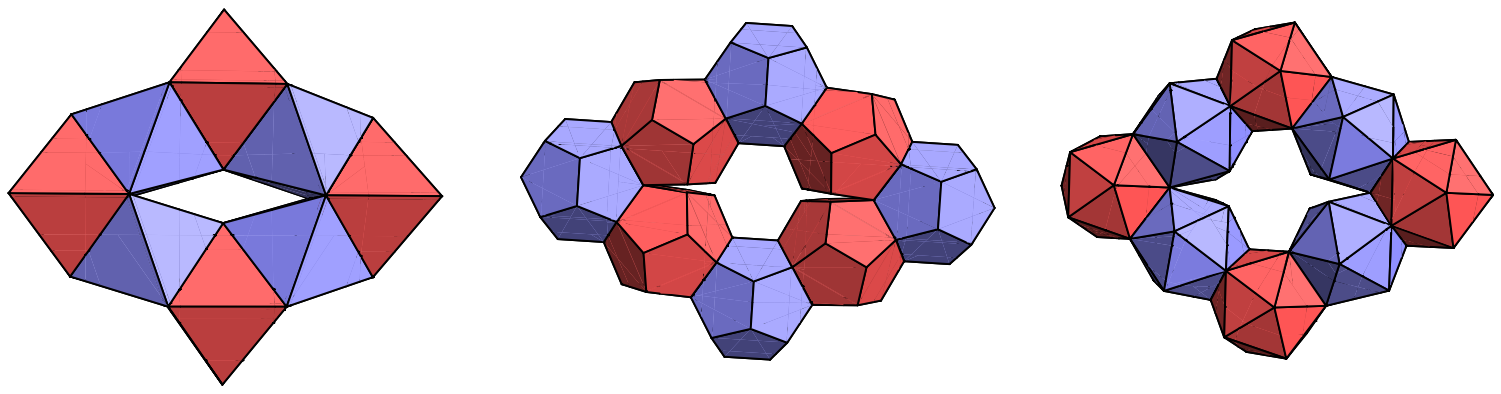

**Figure 1.** Octahedral, dodecahedral, and icosahedral tori.

**Definition 1.** A *tetrahedral chain* is a sequence of $k$ congruent regular tetrahedra meeting face to face, but never doubling back (i.e., $T\,T'\,T$ never occurs). Fixing an ordering of the four vertices of the initial tetrahedron, the chain corresponds to a sequence from $\{1, 2, 3, 4\}$ of length $k-1$, where the integer $i$ denotes a reflection in the $i$th face (the face opposite vertex $i$), and consecutive integers are distinct.

**Theorem 1 (S. Świerczkowski).** The last tetrahedron in a chain $T_0, \ldots, T_n$ cannot coincide with the first.

**Proof.** If $T$ is a regular tetrahedron with vertices $V_i$, let $\phi_i$ be the reflection in the face opposite $V_i$. Any point in $\mathbb{R}^3$ is uniquely representable as $x_1 V_1 + x_2 V_2 + x_3 V_3 + x_4 V_4$, where $\sum x_i = 1$; these are *barycentric coordinates* with respect to $T$. Each $\phi_i$ may be represented by a $4\times4$ matrix $M_i$ acting on these coordinates, where the columns are the vectors $\phi_i(V_j)$; $V M_i$, where $V$ is the $3\times4$ matrix having the $V_j$ as columns, gives the reflected tetrahedron and composition corresponds to matrix multiplication. The reflection in the face $x_i = 0$ sends $V_i$ to $C + (C - V_i) = 2C - V_i = 2\left(\frac{1}{3}\sum_{j\neq i} V_j\right) - V_i$, where $C$ is the centroid of the face opposite $V_i$ (Fig. 2), and so the barycentric matrices for the $\phi_i$ are

$$M_1 = \begin{pmatrix} -1 & 0 & 0 & 0 \\ \frac{2}{3} & 1 & 0 & 0 \\ \frac{2}{3} & 0 & 1 & 0 \\ \frac{2}{3} & 0 & 0 & 1 \end{pmatrix} \quad M_2 = \begin{pmatrix} 1 & \frac{2}{3} & 0 & 0 \\ 0 & -1 & 0 & 0 \\ 0 & \frac{2}{3} & 1 & 0 \\ 0 & \frac{2}{3} & 0 & 1 \end{pmatrix} \quad M_3 = \begin{pmatrix} 1 & 0 & \frac{2}{3} & 0 \\ 0 & 1 & \frac{2}{3} & 0 \\ 0 & 0 & -1 & 0 \\ 0 & 0 & \frac{2}{3} & 1 \end{pmatrix} \quad M_4 = \begin{pmatrix} 1 & 0 & 0 & \frac{2}{3} \\ 0 & 1 & 0 & \frac{2}{3} \\ 0 & 0 & 1 & \frac{2}{3} \\ 0 & 0 & 0 & -1 \end{pmatrix}$$

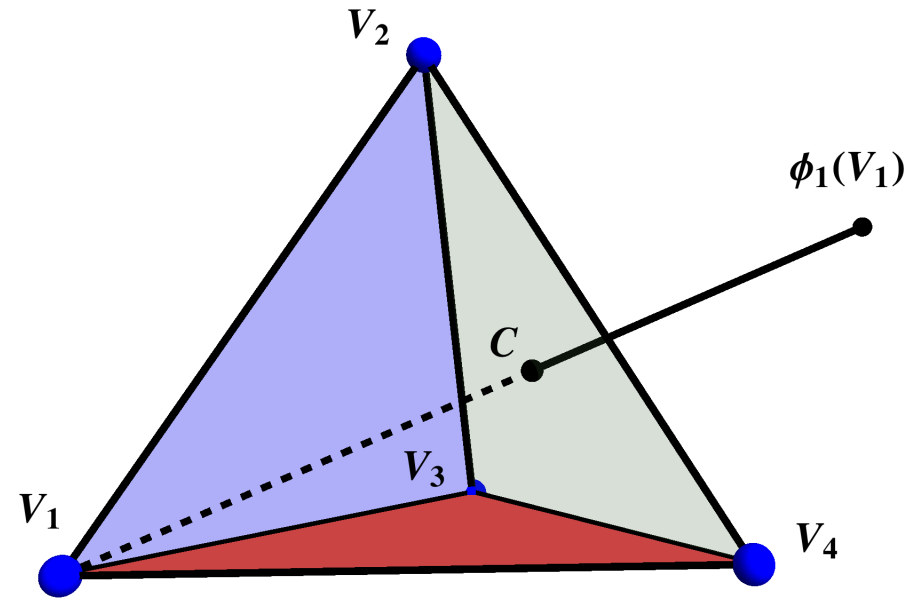

**Figure 2.** Reflection in a side of a regular tetrahedron.

Let $r_0, \ldots, r_{n-1}$ be the reflection sequence of the chain. We may assume $r_0 = 1$. If $T_n$ coincided with $T_0$, then $M_1 M_{r_1} \cdots M_{r_{n-1}}$ is a permutation matrix. The next claim shows that the structure of this matrix product forbids this.

*Claim.* Consider the product $M_1\, M_{r_1} \cdots M_{r_{n-1}}$ with $\frac{2}{3}$ replaced by $x$. The polynomials in the second row have $x$-degree less than $n$ except for the one in the $r_{n-1}$th column, which has degree $n$. And they all have leading coefficient +1.

*Proof of Claim.* By induction on $n$; it is clear for $n = 1$. Consider what happens when the matrix of a word that ends in $M_j$, assumed to have the claimed form, is multiplied on the right by $M_s$, with $s \neq j$. The multiplications by $x$ preserve the claimed property, as the degree becomes $s + 1$ in the $s$th position of row 2, but does not rise at all elsewhere in the row. And the leading coefficient's sign is affected only by the $x$ multipliers.

Now look at the polynomial in the $(2, r_n)$ position: $x^n + a_1 x^{n-1} + \ldots + a_n$, where $a_i \in \mathbb{Z}$. Setting $x = \frac{2}{3}$ and taking a common denominator yields $\left(2^n + 3 a_1 2^{n-1} + \ldots + 3^{n-1} a_n\right) / 3^n$, the numerator of which is not divisible by 3; the fraction is therefore not 0 or 1, as required. □

Perfection may be impossible, but searching for near-perfection is an interesting challenge. In his memoir ([15, p. 191]; see also [16, Conj. 6.1]), Świerczkowski put it this way:

> "Granted then, that the last pyramid in a Steinhaus chain never can have a sidewall in common with the first pyramid $P$, it still may happen that all observations and measurements indicate that these two pyramids do have a sidewall in common. This would not contradict the mathematical result; it would only illustrate the obvious fact that no measurement is 100% accurate. So, a new problem is born: Whatever threshold of accuracy is selected, say, represented by a (small) number $\epsilon$, will there be a chain of pyramids, returning to $P$ such that within the accuracy of $\epsilon$ inches, the last pyramid of the chain has indeed a sidewall in common with $P$. It is hard to tell if anyone will ever want to devote her or his time to search for an answer to this question. In any case, it is unlikely that an answer would be easily found."

In 2015 we showed [6] how various computer searches led to chains with very small error; our best example had 540 tetrahedra with a discrepancy from closure of about $10^{-17}$ (see Fig. 16). That paper posed two challenges: Prove that the error (the discrepancy from a perfect loop) in an embedded tetrahedral chain can be made arbitrarily small; and exhibit specific examples of tetrahedral chains with error under $10^{-18}$. A reexamination of some patterns from [6] led us to the quadrahelix of §3. That chain achieves both goals, by being embedded and having arbitrarily small deviation from exact closure. And despite Świerczkowski's prediction that this resolution will not be easy, the complete proof, once the appropriate pattern has been found, is not very complicated.

The correspondence between chains and strings from {1, 2, 3, 4} needs a little more attention in order to study the gap in a chain. For any chain $C = (T_1, T_2, \ldots, T_n)$ there is a corresponding reflection sequence $r_1, r_2, \ldots, r_{n-1}$. We wish to associate a barycentric matrix to $C$ (as in the proof of Thm. 1), so that the matrix's distance from $I$ is related to the gap in $C$. To this end, we need an invisible starting tetrahedron $T_0$ so that we can study the matrix that takes $T_0$ to $T_n$; if that were the identity then $T_n$ would coincide with $T_0$ and $C$ would close up exactly. So consider some fixed tetrahedron $T_0$ and build the chain by starting with some legal reflection in a face of $T_0$ to get $T_1$ (see Fig. 3). There are three choices because if $r_1$ is, say, 3 then $T_0$ cannot reflect in face 3 to start the chain but it could reflect in faces 1, 2, or 4. Each of these three choices of $r_0$ yields a different reflection string $r_0, r_1, \ldots, r_{n-1}$ and a different placement of the chain in space; in particular, $T_n$ is in a different location for each choice. Each choice yields a barycentric matrix $K$ as the product $M_{r_0} M_{r_1} \ldots M_{r_{n-1}}$. To get $T_n$ from $K$, observe that its vertices are given by $T_0 K$, where, abusing notation slightly, $T_0$ is the $3 \times 4$ matrix whose columns are the vertices of tetrahedron $T_0$.

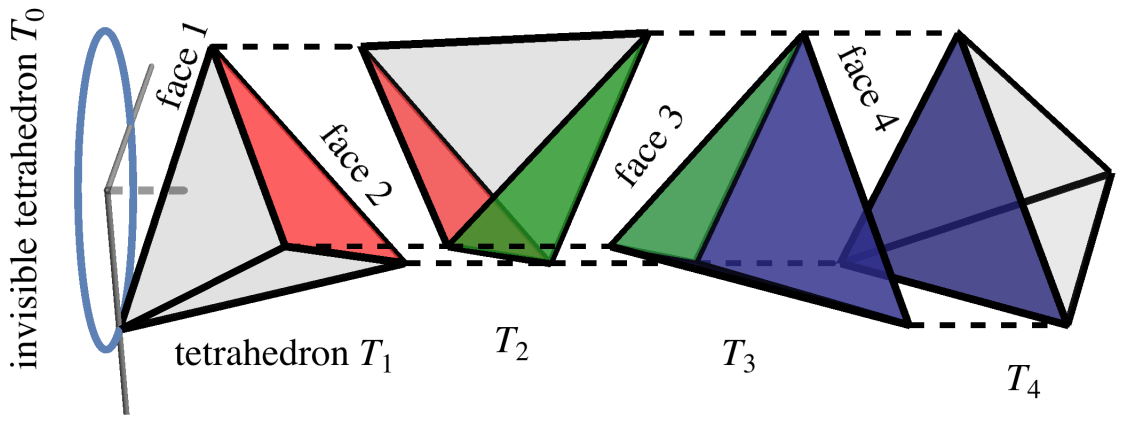

**Figure 3.** The sequence 234 corresponds to a chain of four tetrahedra. There are three choices, 1234, 3234, 4234 for introducing the invisible tetrahedron and the first reflection; the 1234 case is shown.

We need a precise definition of the gap of a tetrahedral chain. There are many ways to do it, all starting from the idea that, for a perfectly closed loop, $T_0$ would coincide exactly with $T_n$. We can measure the difference between these two tetrahedra by looking at the Hausdorff distance $d_H$ between them as sets. For two compact sets $X$ and $Y$, $d_H(X, Y)$ is the smallest $c$ such that an expansion of $X$ by $c$ contains $Y$ and an expansion of $Y$ by $c$ contains $X$ (a $c$-expansion uses radius-$c$ balls around each point of the set). This measure can also be formulated in terms of the distance from a point to a set (see proof of Lemma 2); because our sets are tetrahedra, it is not hard, using some symbolic equation-solving on a derivative, to develop a fast algorithm to compute $d_H$ exactly for two tetrahedra.

An adequate upper bound on $d_H$ is the easily computed $d_H^{\text{discrete}}$, where only the two vertex sets of the tetrahedra are used (a proof uses the alternative form of $d_H$ given in the proof of Lemma 2(a)). A chain has three associated barycentric matrices depending on the choice of initial matrix. The chain closes up perfectly iff one of the three choices of $K$ is a permutation matrix $P$. Thus the minimum of the $3 \cdot 24 = 72$ values of $\|K - P\|_2$ is also an estimate of the discrepancy from perfect closure, where the norm is the standard induced 2-norm; recall that $\|A\|_2$ (often abbreviated to $\|A\|$) is $\sup \|Ax\|/\|x\|$ over nonzero vectors $x$ and equals the largest singular value: the square root of the largest eigenvalue of $A^{\mathrm{T}}A$. Nonidentity permutations do not arise in our work, so we will omit them from the following definition. We will also use $\|K - I\|_{\max}$, the maximum absolute value of the matrix entries.

**Definition 2.** The *gap* of any chain of tetrahedra is the minimum of the three Hausdorff distances between $T_0$ and $T_n$, over the three choices of $r_0$ that lead to three positions for $T_n$. The *norm gap* is $\|K - I\|_2$, again minimized over the three choices of $K$.

The various notions of gap are closely related. We will use matrices so need the following bounds, which relate the gap to the norm gap and to $\|K - I\|_{\max}$. We need a specific invisible tetrahedron, so we will henceforth take $T_0$ to be $\{V_{-1}, V_0, V_1, V_2\}$, where $V_i = \left(\frac{3}{10}\sqrt{3}\cos(i\,\theta),\ \frac{3}{10}\sqrt{3}\sin(i\,\theta),\ \frac{1}{\sqrt{10}}\,i\right)$ with $\theta = \cos^{-1}\left(\frac{2}{3}\right)$ as in the next section.

**Lemma 2.** Consider any tetrahedral chain with invisible tetrahedron $T_0$. Then for any choice of the first reflection $r_0$, with associated barycentric matrix $K$, we have:

**(a).** the gap is no greater than $\|K - I\|_2$;

**(b).** the gap is no greater than $4\|K - I\|_{\max}$.

**Proof. (a)** Let $T_0$ refer to the $3 \times 4$ matrix with the vertices of $T_0$ as columns; let $T_n$ be similar for the final tetrahedron, defined according to the choice of $r_0$. The largest singular value of $T_0$ is $\eta = \sqrt{17}/5$, which is $\|T_0\|_2$. We have $T_0 K = T_n$. Let $\Gamma = \left\{w \in [0, 1]^4 : \sum w_i = 1\right\}$ be the set of all possible barycentric coordinates of points in a tetrahedron. Any point of a solid tetrahedron can be written as the product of that tetrahedron's $3\times 4$ vertex matrix with a vector in $\Gamma$. The Hausdorff distance is

$$d_H(T_0, T_n) = \max[\max_{u\in\Gamma} \min_{w\in\Gamma} \|T_n\, u - T_0\, w\|,\ \max_{w\in\Gamma} \min_{u\in\Gamma} \|T_n\, u - T_0\, w\|];$$

the first term is bounded as follows, where all maxima and minima are over $\Gamma$. A step at the end uses the fact that any barycentric vector lies inside the unit disk.

$$\begin{aligned}\max_u \min_w \|T_n\, u - T_0\, w\| &= \max_u \min_w \|T_0\, K\, u - T_0\, w\|\\ &= \max_u \min_w \|T_0\,(K\, u - w)\| \le \|T_0\| \max_u \min_w \|K\, u - w\|\\ &= \eta \max_u \min_w \|(K-I)\, u - (w-u)\| \le \eta \max_u \min_w \|(K-I)\, u\| + \|w-u\|\\ &\le \eta \max_u \|(K-I)\, u\| \le \eta\, \|K-I\| \max_u \|u\| \le \eta\, \|K-I\| \le \|K-I\|\end{aligned}$$

The bound on the second term is the same, proving (a).

**(b)** For any $n\times n$ matrix $A$, we have the following, where $x$ represents a unit vector and all sums are from 1 to $n$. We use the easily proved fact that $\sqrt{n}$ is the maximum value of $\Sigma_j \left|x_j\right|$ for a unit vector $x$.

$$\begin{aligned}\|A\|^2 &= \max_x \|A\, x\|^2 = \max_x \Sigma_i \left(\Sigma_j\, a_{ij}\, x_j\right)^2 \le \max_x \Sigma_i \left(\Sigma_j \left|a_{ij}\, x_j\right|\right)^2\\ &\le \left(\max_{i,j} \left|a_{ij}\right|^2\right) \max_x \Sigma_i \left(\Sigma_j \left|x_j\right|\right)^2 \le \|A\|_{\max}^2\, \Sigma_i \left(\sqrt{n}\right)^2 = n^2\, \|A\|_{\max}^2\end{aligned}$$

When $n=4$, this is $\|A\| \le 4\, \|A\|_{\max}$, which, by (a), gives (b). □

Many of our proofs rely on algebraic computation and manipulation of determinants, normal vectors, and so on. The relevant *Mathematica* code and some intermediate formulas are in an Appendix. For a slightly different approach to the main result and its proof, using more geometry and less algebra, see [7].

**Acknowledgment.** We thank Danny Lichtblau of Wolfram Research, Inc., and Damien Stehlé for providing the details of a lattice-reduction solution to a Diophantine approximation problem. And we are indebted to the graphic and algebraic capabilities of *Mathematica*, without which our discoveries would not have been possible.

## 2. The Boerdijk–Coxeter Helix

Our basic building block is the tetrahelix (also called the Boerdijk–Coxeter helix), a stack of regular tetrahedra with some pleasant properties [2, 8]. Let $r = \frac{3}{10}\sqrt{3}$, $h = \frac{1}{\sqrt{10}}$, and $\theta = \cos^{-1}\left(-\frac{2}{3}\right) \sim 132\,°$. The *tetrahelix* of length $L$ is the linear chain of $L$ tetrahedra defined by the $L+3$ points $\{V_i : i = 0, 1, \ldots, L+2\}$, given in cylindrical form as $V_i = (r\cos(i\,\theta), r\sin(i\,\theta), i\,h)$; Figure 4 shows how the vertices spiral up along the radius-$r$ helix centered on the $z$-axis. Each tetrahedron has side-length 1. The reflection sequence for the tetrahelix is the periodic form 1234 1234 1234 12.…

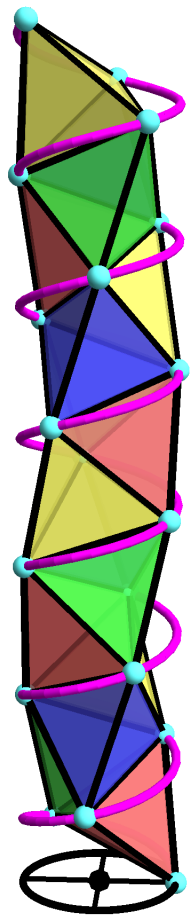

**Figure 4.** The tetrahelix made from 16 tetrahedra colored, in order, red, green, blue, yellow, red, green, blue, yellow, and so on. The vertices are equally spaced along a helix.

Each point $V_q$ can be represented barycentrically in terms of the basic points $T_0 = \{V_{-1}, V_0, V_1, V_2\}$ (the invisible tetrahedron as in §1). Using the cylindrical formula, one can compute these barycentric coordinates explicitly. Then in Cartesian coordinates $V_q = T_0\, C$, where $T_0$ is viewed as a $3 \times 4$ matrix with columns $V_i$ $(-1 \le i \le 2)$ and $C$ is the $4 \times 1$ matrix given by

$$(*) \qquad C = \frac{1}{10}\begin{pmatrix} 3 \\ 4 \\ 3 \\ 0 \end{pmatrix} + \frac{1}{10}\begin{pmatrix} -3 \\ -1 \\ 1 \\ 3 \end{pmatrix} q + \frac{3}{10}\begin{pmatrix} -1 \\ 2 \\ -1 \\ 0 \end{pmatrix} \cos(q\,\theta) + \frac{3\sqrt{5}}{50}\begin{pmatrix} -2 \\ 1 \\ 4 \\ -3 \end{pmatrix} \sin(q\,\theta)$$

Note that each of the last three column vectors sums to 0 while the first sums to 1, and hence $C$ sums to 1, as is always the case for barycentric coordinates. The tetrahelix has been used in one unusual building construction; in 1989 Isozaki Arata designed a 100-meter high tower (Fig. 5) in Mito, Japan, in the shape of the tetrahelix with 28 tetrahedra [1].

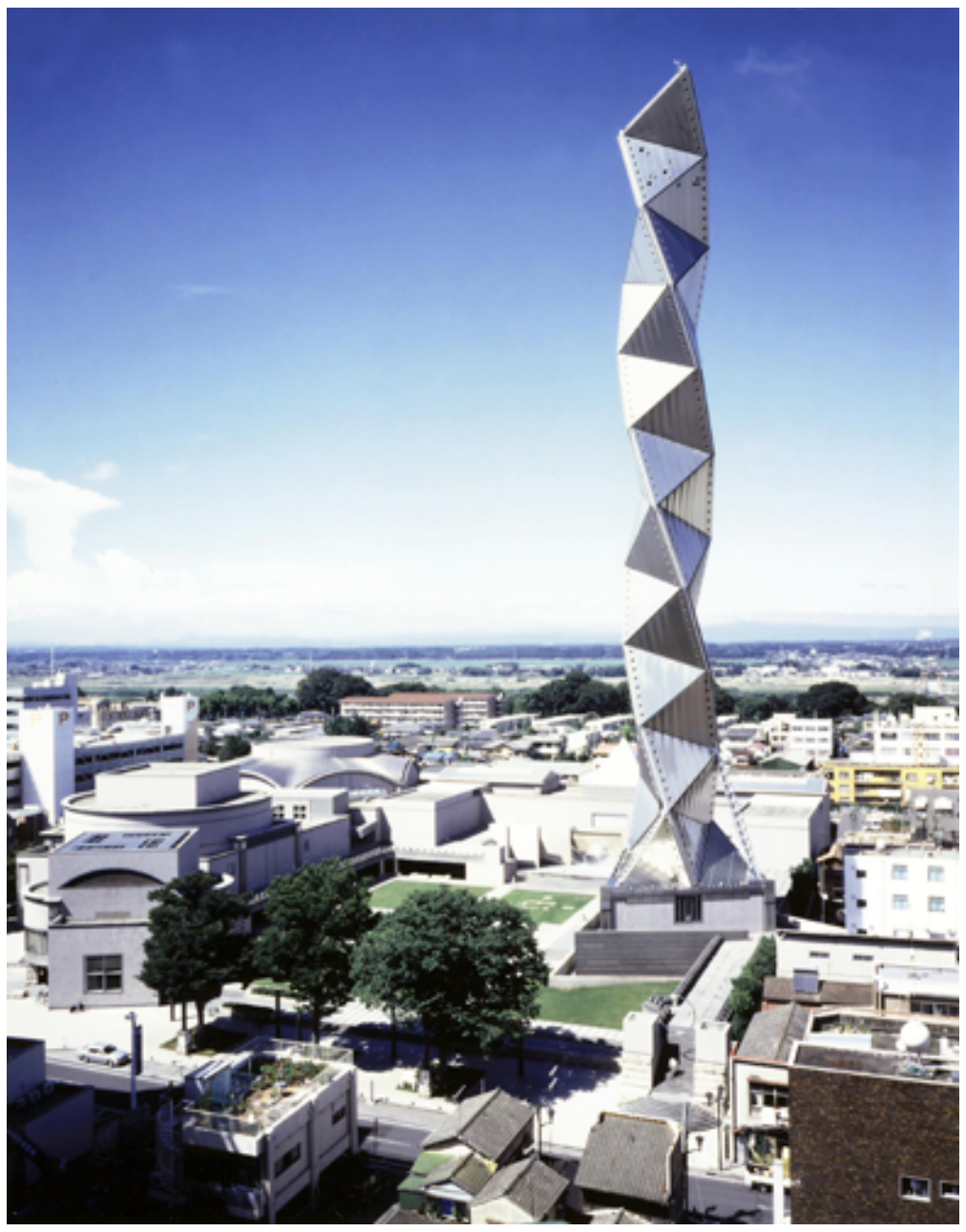

**Figure 5.** A giant tetrahelix rises from the art museum in Mito, Japan.

## 3. The Quadrahelix

The *quadrahelix* $QH_L$ is built by linking four tetrahelices of length $L + 1$ using a common tetrahedron at the first and third joins and face-attachment at the second; $QH_L$ has $4L + 2$ tetrahedra. To get the appropriate reflection string, let $S_m$ denote the $m$-term tetrahelix string that begins with a 2: 2341234.... Let $\Sigma_m$, for even $m$, be $S_{m+1}$ with its middle entry deleted; such a deletion is a way of making the appropriate first turn. Then $QH_L$ is simply $1\, \Sigma_{2L}\, j\, \overline{\Sigma_{2L}}$ where $j = 3$ for even $L$ and 1 for odd $L$ and $\bar{s}$ is the reversal of $s$. The corresponding barycentric matrix is the product $M_1\, M_2\, M_3 \cdots M_3\, M_2$ defined from the full string. Thus $QH_4$ is 1 2341 3412 3 2143 1432 and $QH_{10}$ is 1 2341234123 1234123412 3 2143214321 3214321432. The deletions have the effect of making $T_{L+1}$ and $T_{3L+2}$ into pivots; they are each part of two tetrahelices. Figure 6 shows the chains for $L = 5, 6, 7, 10$; the two pivot tetrahedra are shown in gray and so each colored section has $L$ tetrahedra.

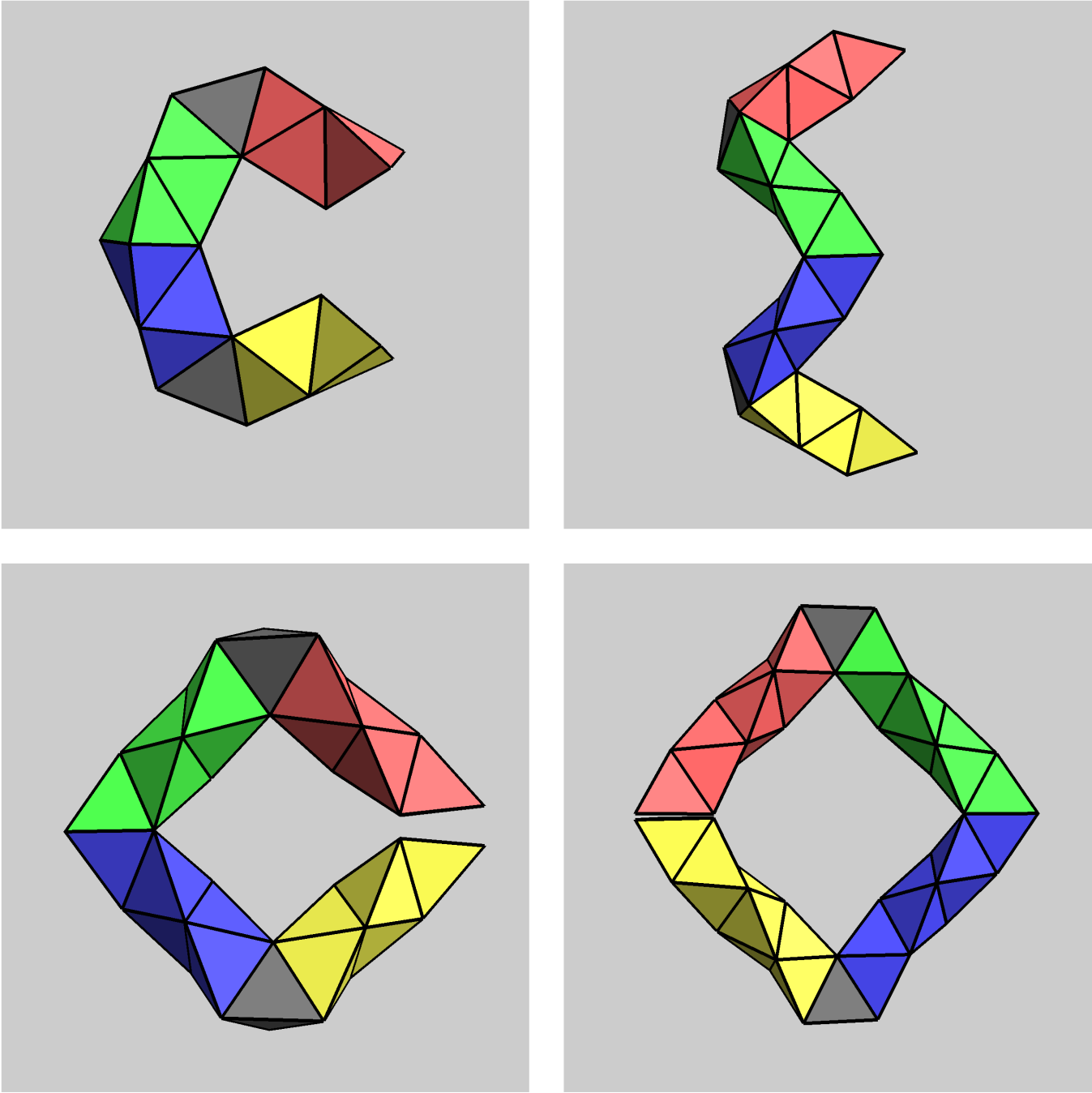

**Figure 6.** $QH_5$, $QH_6$, $QH_7$, and $QH_{10}$. The gap for $QH_{10}$ is 0.078, about 8% of a tetrahedon side. There are no collisions in these, and indeed none in any quadrahelix.

The basic tetrahelix is a surprising shape, looking essentially the same regardless of its length. Remarkably, the quadrahelix shares this property, in the sense that for any $L$, the quadrahelix forms a 4-sided path having no collisions; the leg lengths (viewed along the tetrahelix axes) are all equal; and the three angles between the tetrahedral axes are $\sec^{-1}(5)$, $\sec^{-1}(-5)$, and $\sec^{-1}(5)$. It exhibits several useful symmetries (see Fig. 10). The main result of this paper, that $QH_L$ can be arbitrarily close to a closed loop, will require a proof of embeddability together with sufficient analysis to guarantee that there are values of $L$ for which the final gap is arbitrarily small. There is in fact a single shape—the limiting rhombus—that the quadrahelices converge to, for special values of $L$ (Fig. 13). Figure 7 shows $QH_{70}$, a typical almost-closed chain; the gap (between red and yellow) is less than 1% of a tetrahedron side and is not visible at the full scale. For $QH_{1960}$, the gap is $\frac{1}{100}$% (of a tetrahedron edge). And when $L$ is the 99-digit integer

$$521\,269\,338\,782\,055\,651\,792\,691\,214\,128\,196\,053\,088\,348\,030\,247\,372\,007\,924\,246\,566\,932\,650\,514\,801\,545\,115\,813\,925\,856\,156\,787\,510,$$

the gap size is under $10^{-101}$.

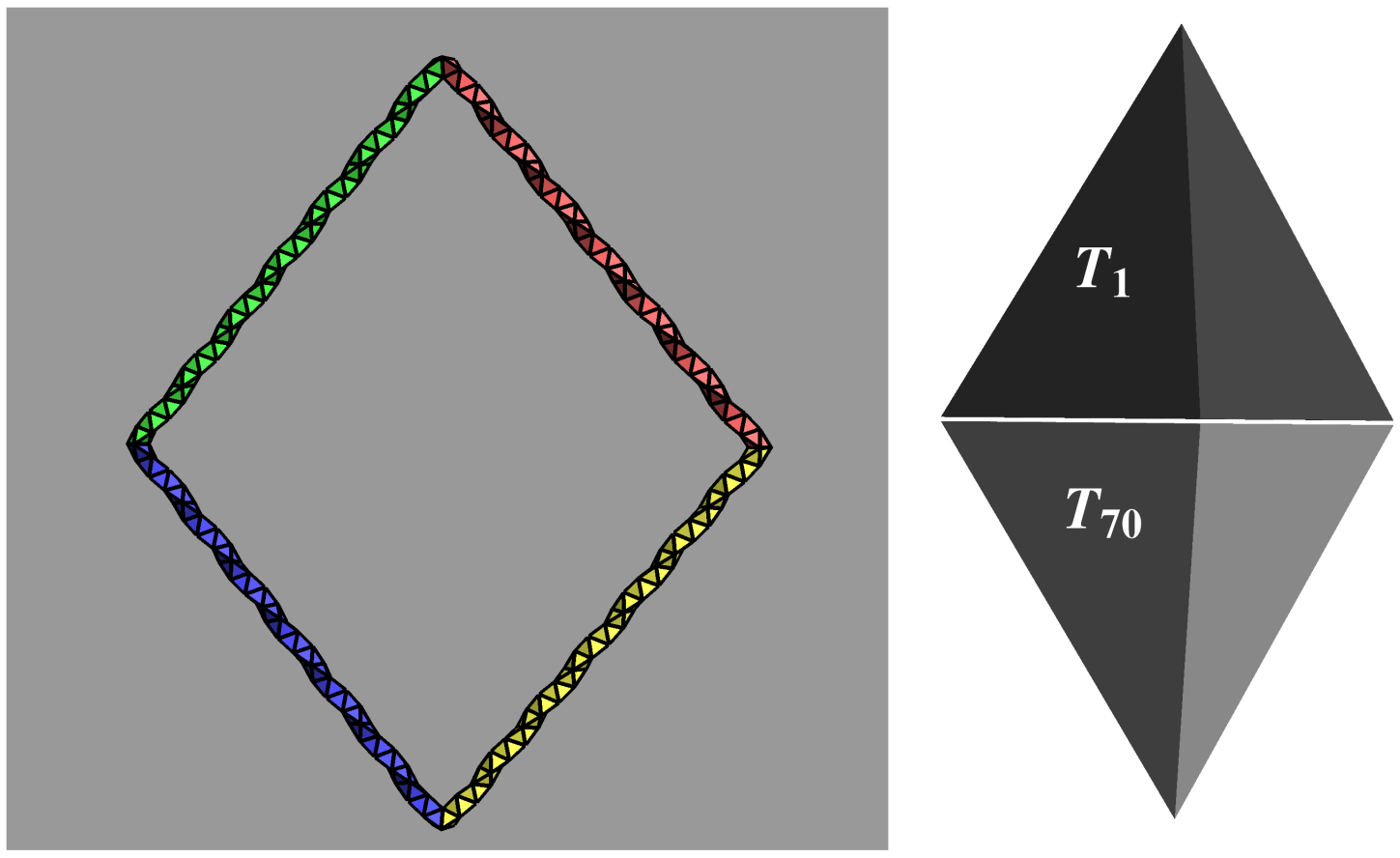

**Figure 7.** $QH_{70}$ has a gap less than 1% of a tetrahedron side. The gap between red and yellow is shown in the magnified image at right.

Our main proof in §4 will use some exact algebraic formulas, but the underlying geometry that makes the quadrahelix work as a loop is easy to understand. Each new vertex of the tetrahelix rotates $\theta$ around the axis. If $(L+1)\,\theta$ is close to a multiple of $2\pi$, then $T_{L+1}$ is close to being a translation of $T_0$: the invisible tetrahedron described in §1. Therefore the starting triangle—the floor of $T_1$—is almost parallel to the ceiling of $T_{L+1}$ (Fig. 8). This in turn means that any plane orthogonal to this ceiling makes almost a right angle with the start triangle. The plane that bisects $T_{L+1}$ as in Figure 10 (the first quadplane) then serves as a reflecting plane for the initial tetrahelix, and two more reflections cause a total change of nearly 360°, making a nearly closed loop.

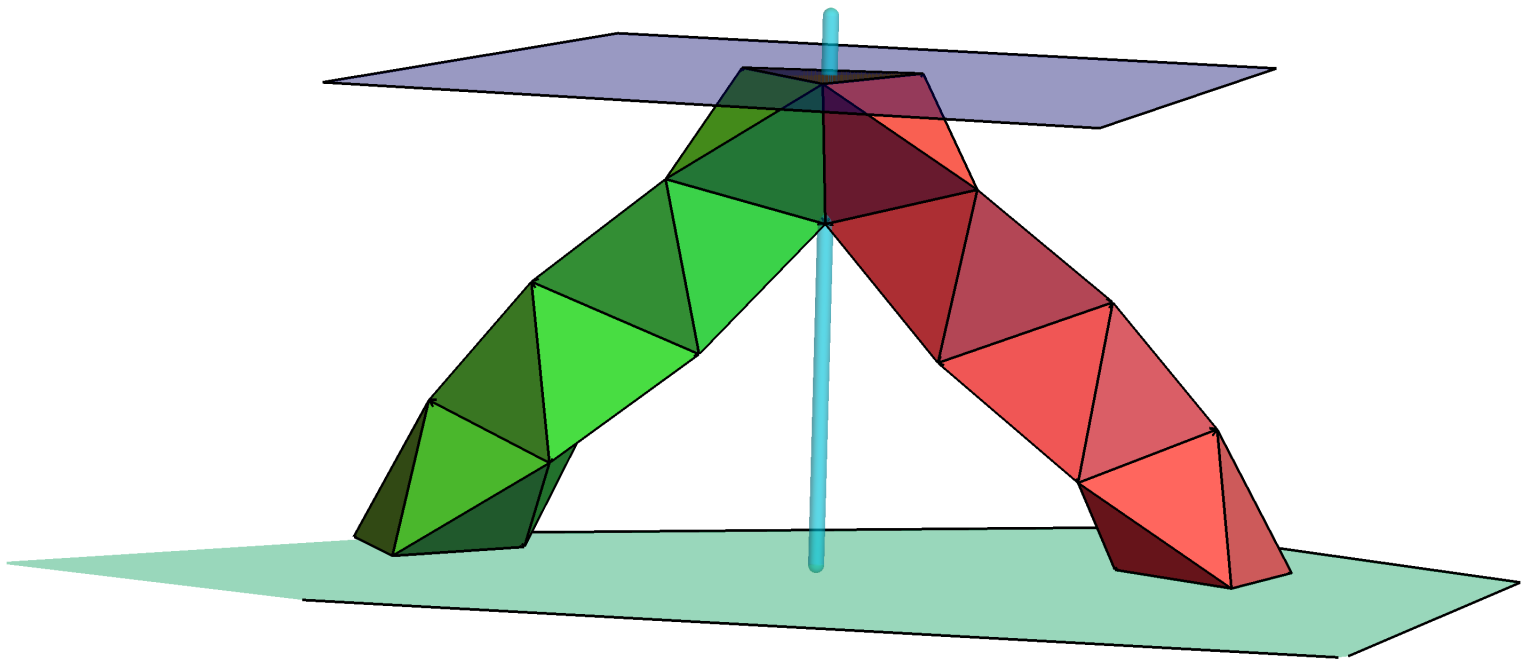

**Figure 8.** For certain $L$ (such as $L = 10$) the tetrahelix has a ceiling (upper plane) that is nearly parallel to the floor (plane through initial red triangle).

One can use 3D printing technology to realize these chains. Figure 9 shows a Shapeways model of $QH_{10}$, with the gap in black; this model is available at [17]. The quadrahelix pattern is so simple that it would seem not difficult to build a gap-free model of $QH_{10}$, $QH_{29}$, or $QH_{40}$ using a standard polyhedron construction tool such as *ZomeTool* or *Polydron*.

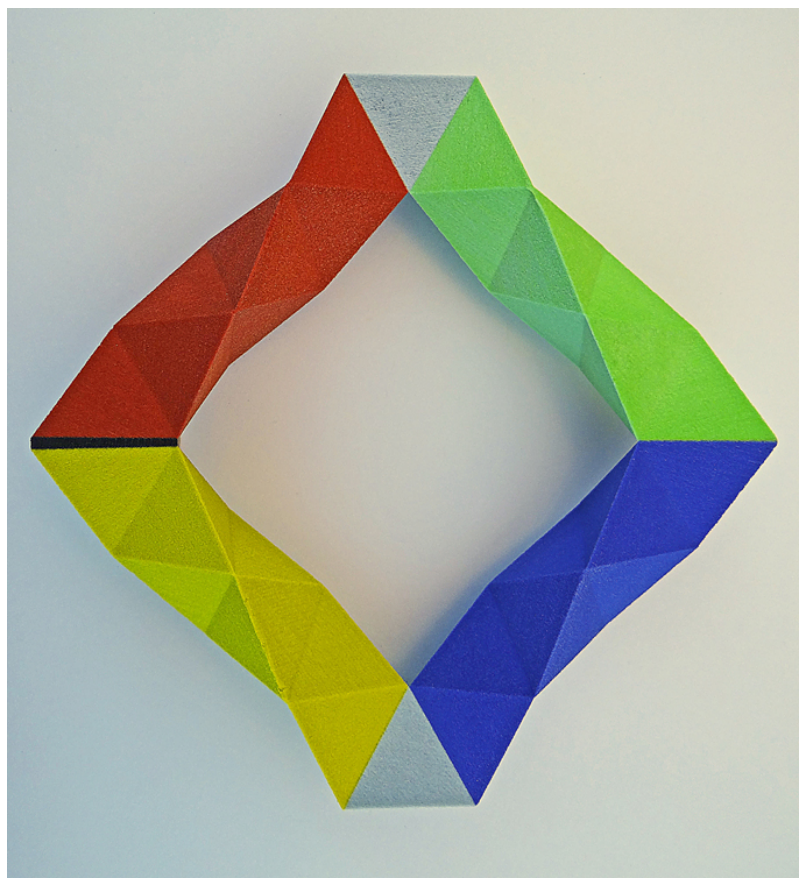

**Figure 9.** A colored sandstone model of $QH_{10}$, printed by Shapeways.

## 4. The Disappearing Quadrahelix Gap

Our main result is that $QH_L$ is always embedded and can have an arbitrarily small gap. To prove the chains are embedded, consider $QH_L$ as falling into four congruent parts, which we think of as red, green, blue, and yellow, with red being the start of the chain. Think of the first pivot tetrahedron (gray in Fig. 10) as being half red and half green, and similarly for the pivot separating blue from yellow. The main dividing plane is the *biplane*, defined by the exact central triangle (the green–blue boundary in Fig. 10). The *first quadplane* is the one that splits the chain that lies on the side of the biplane nearest the start into two equal parts: thus it splits $T_{L+1}$ exactly in half. The *second quadplane* is similar, on the biplane's other side. Thus each colored sector has $L + \frac{1}{2}$ tetrahedra. The reflection string for $QH_L$ is highly symmetrical and that combinatorial symmetry corresponds to geometric symmetry in the chain.

**Reflection Lemma.** The second tetrahelix of $OH_L$ is a reflection of the first in the first quadplane (and similarly for the second quadplane). The second half of the chain is a reflection of the first in the biplane.

**Proof.** The reflection faces used to enter and exit the pivot tetrahedron are $L+1$ and $L+3$ (such indices are always considered as reduced mod 4). Therefore the reflection string leaving the pivot in the start direction is $L+1, L+4, L+3, L+2, L+1, \ldots$ and in the other direction is $L+3, L, L+1, L+2, L+3, \ldots$. The second is the result of applying an $L+1 \leftrightarrow L+3$ switch to the first. But such a switch corresponds to a reflection in the quadplane, since such a reflection exchanges faces $L+1$ and $L+3$ at the pivot, and hence throughout. The biplane assertion is a simple consequence of the palindromic property of the string defining $QH_L$. The pivot object in this case is the face defined by the central entry of the string. □

**Five-Plane Lemma.** The biplane, first quadplane, second quadplane, start plane, and end plane all pass through a common line, called the *rotation axis*.

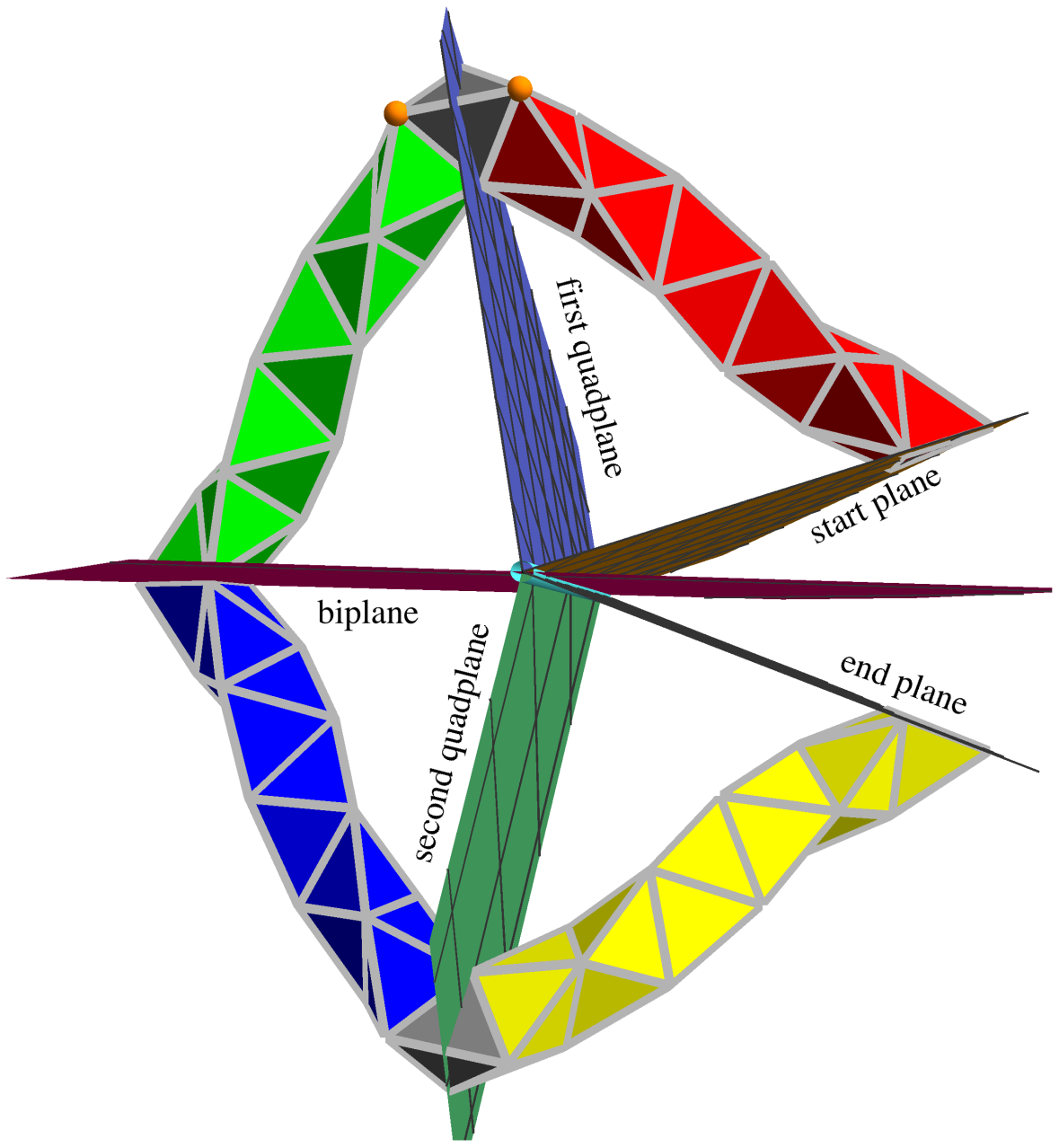

**Figure 10.** The quadrahelix has three symmetry planes. These planes, and the planes defined by the two ends, all pass through a common rotation axis. Shown here is $QH_{15}$; the two orange points at top are $V_{L+1}$ and $V_{L+3}$; the vector they determine is normal to the first quadplane.

**Proof.** The start plane meets the first quadplane in a line; by the Reflection Lemma, the start plane and the biplane are reflections in the first quadplane; therefore the biplane meets the first quadplane in the same line. Reflection in the biplane yields the remaining collinearities. □

Let $\rho_0$ be the dihedral angle between the start plane and the first quadplane, in the sector containing the start of the quadrahelix; let $\rho = \frac{\pi}{2} - \rho_0$ (see Fig. 11).

**Acute Angle Lemma.** Angle $\rho_0$ is acute.

The first quadplane is defined by the bisected tetrahedron atop the first tetrahelix. This means its normal vector is the edge of that tetrahedron connecting the tetrahelix vertices $V_{L+1}$ and $V_{L+3}$ (see Fig. 10). Thus the preceding lemma is equivalent to the next one. We define $\overline{\alpha}$ to be the unique value in $[-\pi, \pi)$ that is congruent to $\alpha$ modulo $2\pi$; and we will use $\delta$ for $(L+1)\,\theta$.

**Acute Angle Lemma.** The angle between $A$, the vector normal to the tetrahelix's base triangle $V_0\, V_1\, V_2$ and pointing in the direction of $V_3$, and the vector $V_{L+3} - V_{L+1}$ is acute.

**Proof.** A purely geometric proof can be found, but algebra is quicker. First, $A = \frac{1}{\sqrt{6}}(3\, V_3 - (V_0 + V_1 + V_2)) = \frac{1}{3\sqrt{5}}\left(\sqrt{10}, 2\sqrt{2}, 3\sqrt{3}\right)$. Using the substitution $L\theta \to \delta - \theta$ and some trig simplification yields $V_{L+3} - V_{L+1} = \frac{-1}{\sqrt{15}}\left(\sqrt{5}\cos\delta - 2\sin\delta, 2\cos\delta + \sqrt{5}\sin\delta, \sqrt{6}\right)$ and $A\cdot(V_{L+3} - V_{L+1}) = \frac{\sqrt{6}}{5}(1 - \cos\delta) = \frac{2\sqrt{6}}{5}\sin^2\left(\frac{1}{2}\delta\right)$. If $\delta$ is an integer multiple of $2\pi$, then the dot product vanishes, the angle in question is exactly 90°, and the loop closes up perfectly, contradicting Theorem 1. Therefore the dot product satisfies $\frac{2}{5}\sqrt{6} \geq A\cdot(V_{L+3} - V_{L+1}) > 0$ and the angle is strictly between 11° and 90°. □

Instead of Theorem 1, we could have used the irrationality of $\frac{\theta}{2\pi}$ [12, Cor. 3.12] to deduce that $\delta$ is not a multiple of $2\pi$. In fact, the use of Theorem 1 yields an alternative proof of the irrationality result. An important corollary of the preceding proof is the following, which shows the direct connection between $\rho$ and the angle corresponding to $(L+1)\,\theta$: one is small if and only if the other is, and the relationship is quadratic. Recall that $\delta$ abbreviates $(L+1)\,\theta$, and an overbar denotes the reduced angle modulo $2\pi$.

**Corollary 3(a).** When $\rho_0$ is the dihedral angle between the biplane and first quadplane and $\rho = \frac{\pi}{2} - \rho_0$, we have $\cos\rho_0 = \sin\rho = \frac{2\sqrt{6}}{5}\sin^2\frac{\delta}{2}$.

**(b).** If $R \in SO_3(\mathbb{R})$ is the 3-dimensional rotation matrix through angle $4\rho_0$, then $\|R - I\| < \bar{\delta}^2$.

**Proof.** (a) follows from the dot product in the acute angle lemma. (b) follows from $4\rho_0 = \cos^{-1}\left(\frac{2}{5}\sqrt{6}\sin^2\frac{\delta}{2}\right)$ and the fact that the norm of the difference between an angle-$\alpha$ rotation matrix and the identity matrix is $2\left|\sin\frac{\alpha}{2}\right|$. □

**Theorem 4.** For any $L$, the quadrahelix $\text{QH}_L$ is embedded.

**Proof.** A tetrahelix has no collision, so each colored sector is embedded. The next step is showing that the first sector (red) stays within the region defined by the start plane and the quad plane. The Reflection Lemma then yields the same for the other colors and the appropriate planes. For this proof, we can view the pivot tetrahedron attached to the red tetrahelix as being the start; then the three points defining the quadplane are $V_1$, $V_2$, and $\frac{1}{2}(V_0 + V_3)$. The sign of $D$, the determinant of the $3\times 3$ matrix $\left(V_q - V_1, V_q - V_2, V_q - \frac{1}{2}(V_0 + V_3)\right)$, determines which side of the quadplane contains $V_q$. Letting $c = \cos(q\,\theta)$ and $s = \sin(q\,\theta)$, the matrix is this:

$$r\begin{pmatrix} c + \frac{2}{3} & s - \frac{\sqrt{5}}{3} & (q-1)\,\frac{h}{r} \\ c + \frac{1}{9} & s + \frac{4\sqrt{5}}{9} & (q-2)\,\frac{h}{r} \\ c - \frac{49}{54} & s - \frac{7\sqrt{5}}{54} & \left(q - \frac{3}{2}\right)\frac{h}{r} \end{pmatrix}$$

and the value of $20\sqrt{10}\,D$ works out to $6\sqrt{5}\,q - 9\sqrt{5} - \sqrt{5}\,c + 7\,s$, which is not less than $13q - 30$. Because $q \geq 3$ in this approach, $D$ is positive, as required. A similar determinant computation shows that the first tetrahelix stays on the proper side of the start plane. The Five-Plane and Acute Angle Lemmas mean that the five planes define four regions containing the four colored sectors of the quadrahelix and there are no collisions between colors. In particular, the total angle as one moves around the rotation axis is $4\rho_0 < 360°$ and the tetrahedra at the ends of the chain are disjoint. □

Theorem 1 shows that the gap cannot be 0. It is not hard to find *nonembedded* chains with arbitrarily small gaps [6, p. 61]. The fact that the quadrahelix is both embedded and achieves vanishingly small gaps is remarkable, considering how simple the chain is.

The next theorem concludes the proof that embedded tetrahedral chains achieve vanishingly small gaps. Our proof is algebraic, but the heart of the matter is really geometric. Figure 11 shows how the size of the jaw that, in essence, defines the gap, is close to $2\,H\sin(2\rho)$, where $H$ is the distance (projected onto the plane perpendicular to the rotation axis) from the rotation axis to the first tetrahedron. When the jaw is small $H \leq (4\,L + 2)\,h$ and, by Corollary 3, $\sin\rho \leq \frac{\sqrt{6}}{10}\bar{\delta}^2$; these facts mean that the jaw size is bounded by $2\,L\,\bar{\delta}^2$, which, as shown in the proof that follows, approaches 0 for an infinite subsequence of $L$-values.

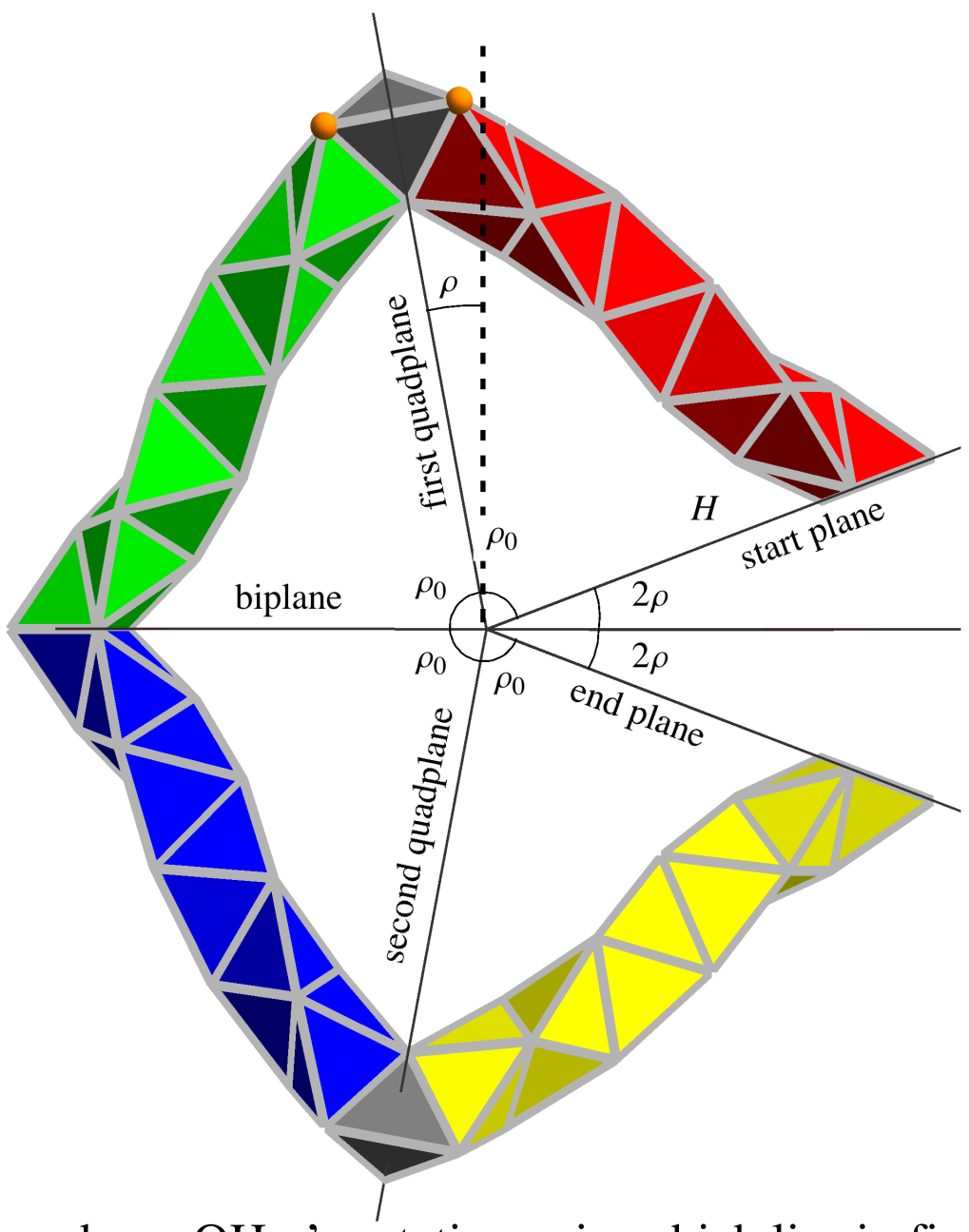

**Figure 11.** Looking down $QH_{15}$'s rotation axis, which lies in five planes. Angle $\rho_0$ is always acute so $\rho$ is positive.

**Theorem 5.** For any $\epsilon > 0$, there is a quadrahelix $QH_L$ having gap less than $\epsilon$.

**Proof.** Recall that $\delta = (L+1)\,\theta$. If $\overline{\delta}$ were 0, then, by Corollary 3, $\rho_0$ would be exactly 90° and $QH_L$ would close up perfectly: the four right angles would form a loop with zero gap. This cannot happen (by either Thm. 1 or the irrationality of $\theta/\pi$), but we can try to make $\overline{\delta}$ small. The convergents $\frac{k}{q}$ of the continued fraction of $\frac{\theta}{2\pi}$ give values $q = L+1$ for which $\frac{(L+1)\,\theta}{2\pi}$ is very close to the integer $k$. Moreover, the error $\left|\frac{\theta}{2\pi} - \frac{k}{q}\right|$ is well known to be bounded by $\frac{1}{\sqrt{5}\,q^2}$ (Hurwitz's Thm.; [4, Exer. 7.10]). Multiplication by $2\pi q$ implies that $\left|\overline{(L+1)\,\theta}\right| < \frac{2\pi}{\sqrt{5}\,(L+1)}$. Using this method, choose $L$ so that $L > \frac{2}{\epsilon}$ and $|\overline{\delta}| < \frac{2\pi}{\sqrt{5}\,(L+1)}$. An alternative to continued fractions is Kronecker's theorem, which gives infinitely many $L$ so that $|\overline{\delta}| < \frac{3}{L}\,2\pi$ [10, Thm. 440].

We can use the barycentric formula $(*)$ of §2 to get an exact symbolic expression for $K$, the matrix giving the barycentric coordinates of the final tetrahedron of $QH_L$. This is done in four steps. The base points are $\text{base}_1 = T_0 = \{V_{-1}, V_0, V_1, V_2\}$, where $V_i$ are from the cylindrical formula for the tetrahelix. The base for the second leg (the tetrahedra after the first turn; green in Fig. 12) is $\text{base}_2 = \{V_{L+4}, V_{L+1}, V_{L+2}, V_{L+3}\}$. These points are used as the base in the barycentric formula and lead to the general point on the second leg, which in turn yields $\text{base}_3$: the points on the second leg corresponding to basic tetrahelix point $\{V_{L+3}, V_{L+2}, V_{L+1}, V_{L+4}\}$. This leads to the third leg (blue), the fourth base (which uses the same permutation as $\text{base}_2$), and the final leg and final tetrahedron.

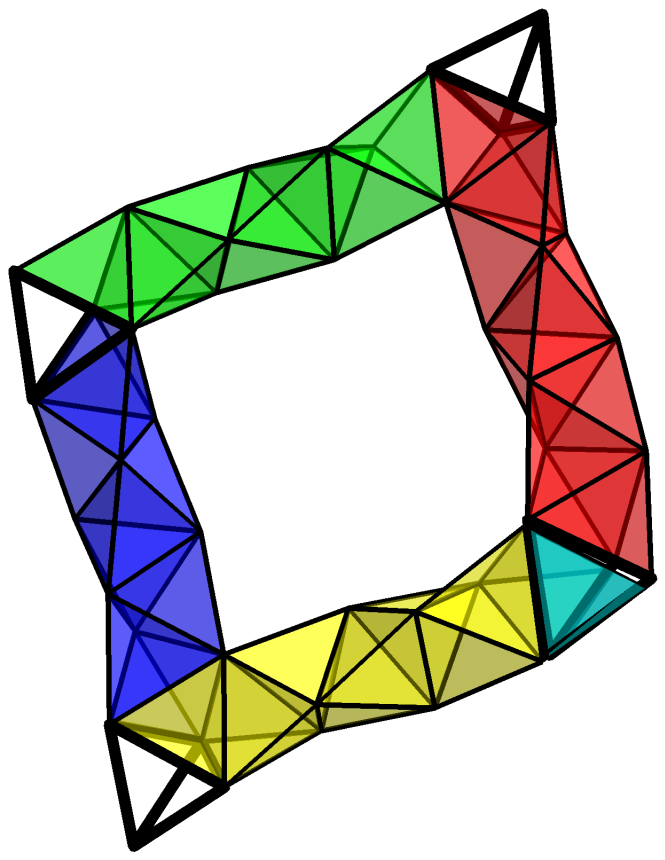

**Figure 12.** $\mathrm{QH}_{10}$, with the four barycentric base sets for the tetrahelices shown as black edges. The one between green and blue ($\mathrm{base}_3$) is in the quadrahelix; the others are not. The cyan tetrahedron is the final one in the chain.

The resulting reasonably concise formula is $K = I + \sigma\,\frac{4}{3125}\left(H_0 + H_1\,\sigma + H_2\,\sigma^2 + H_3\,\sigma^3\right)$, where $\sigma = \sin^2\!\left(\frac{\delta}{2}\right)$ and

$$H_0 = 25\left(5\begin{pmatrix} 6 & 21 & 6 & -9 \\ -22 & -7 & -2 & 3 \\ -2 & -7 & -2 & 3 \\ 18 & -7 & -2 & 3 \end{pmatrix} + 10\begin{pmatrix} 3 & 3 & 3 & 3 \\ -1 & -1 & -1 & -1 \\ -1 & -1 & -1 & -1 \\ -1 & -1 & -1 & -1 \end{pmatrix} L + 3\sqrt{5}\begin{pmatrix} -1 & -6 & 9 & -6 \\ 7 & 2 & -3 & 2 \\ -13 & 2 & -3 & 2 \\ 7 & 2 & -3 & 2 \end{pmatrix}\sin\delta\right)$$

$H_1 =$

$$2\left(25\begin{pmatrix} -120 & -45 & 0 & 45 \\ 76 & -69 & -24 & 21 \\ 40 & 15 & 0 & -15 \\ 4 & 99 & 24 & -51 \end{pmatrix} + 600\begin{pmatrix} 0 & 0 & 0 & 0 \\ -1 & -1 & -1 & -1 \\ 0 & 0 & 0 & 0 \\ 1 & 1 & 1 & 1 \end{pmatrix} L + 3\sqrt{5}\begin{pmatrix} 0 & 0 & 0 & 0 \\ 6 & 26 & -14 & 6 \\ -8 & -28 & -8 & 12 \\ 2 & 2 & 22 & -18 \end{pmatrix}\sin\delta + 120\sqrt{5}\begin{pmatrix} 0 & 0 & 0 & 0 \\ 1 & 1 & 1 & 1 \\ -2 & -2 & -2 & -2 \\ 1 & 1 & 1 & 1 \end{pmatrix} L\sin\delta\right)$$

$$H_2 = 48\left(5\begin{pmatrix} -6 & -21 & -6 & 9 \\ 26 & 31 & 26 & -39 \\ 6 & 31 & -34 & 21 \\ -26 & -41 & 14 & 9 \end{pmatrix} + 10\begin{pmatrix} -3 & -3 & -3 & -3 \\ 4 & 4 & 4 & 4 \\ 1 & 1 & 1 & 1 \\ -2 & -2 & -2 & -2 \end{pmatrix} L + 3\sqrt{5}\begin{pmatrix} 1 & 6 & -9 & 6 \\ -8 & -13 & 12 & -3 \\ 13 & 8 & 3 & -12 \\ -6 & -1 & -6 & 9 \end{pmatrix}\sin\delta\right)$$

$$H_3 = 1440\begin{pmatrix} 4 & 3 & 0 & -3 \\ -5 & -2 & -3 & 6 \\ -2 & -5 & 6 & -3 \\ 3 & 4 & -3 & 0 \end{pmatrix}$$

Using this formula and assuming $L > \frac{2}{\epsilon}$, $|\bar{\delta}| < \frac{2\pi}{\sqrt{5}\,(L+1)}$, and also $\epsilon < \frac{1}{100}$, we next show that $4\,\|K - I\|_{\max} \le \epsilon$. We can bound $\sin\delta$ by $|\bar{\delta}|$ and $\sin^2\!\left(\frac{\delta}{2}\right)$ by $\frac{1}{4}\,\bar{\delta}^2$; then the assumptions on $L$ and $\delta$ yield bounds on each of the 16 entries in $K - I$ that are polynomial in $\epsilon$. These polynomials are easily bounded in absolute value by replacing negative coefficients by positive. For example, the $(1, 1)$th entry, after replacing each symbolic coefficient (they involve rationals, $\sqrt{5}$, and powers of $\pi$) by a slightly larger approximate real is bounded by

$$\frac{\epsilon}{4}\left(0.95 + 3.43\,\epsilon + 5.3\,\epsilon^{2} + 6.42\,\epsilon^3 + 6.53\,\epsilon^4 + 4.25\,\epsilon^5 + 1.45\,\epsilon^6 + 0.64\,\epsilon^7\right).$$

The degree-7 polynomial is under 1 when $\epsilon < \frac{1}{100}$ giving the desired bound of $\frac{\epsilon}{4}$. More detail is in the Appendix. Then Lemma 2(b) asserts that the gap of $\mathrm{QH}_L$ is bounded by $4\,\|K - I\|_{\max} < \epsilon$, as required. □

The proof yields more information. First, the barycentric formula leads to the ideal rhombus that the chains approach as $L$ increases through values so that $\overline{(L+1)\,\theta} \to 0$. We can choose points $P_1, P_2, P_3$ on the first three bases, do some trig expansion, replace $(L+1)\,\theta$ with 0 and normalize the differences $P_2 - P_1$ and $P_3 - P_2$ to get $\frac{1}{\sqrt{(L+4)^2+1}}\left(-\frac{\sqrt{5}}{3\sqrt{6}}, \frac{7}{3\sqrt{6}}, L+4\right)$ and $\frac{-1}{\sqrt{L\,(L+8)+17}}\left(\frac{12\,L+37}{3\sqrt{30}}, \frac{24\,L+119}{15\sqrt{6}}, \frac{L+6}{5}\right)$. The dot product of these is $-\frac{1}{5} - \frac{19+4\,L}{85+40\,L+5\,L^2}$, which approaches $-\frac{1}{5}$ as $L \to \infty$. It follows that the four angles of the limiting rhombus are $\alpha^+, \alpha^-, \alpha^+, \alpha^-$ where $\alpha^{\pm} = \sec^{-1}(\pm 5)$ and it is then routine to find that the normalized rhombus (Fig. 13) has vertices $(0, 0)$, $(0, 1)$, $\frac{1}{5}\left(2\sqrt{6}, 4\right)$, $\frac{1}{5}\left(2\sqrt{6}, -1\right)$. A consequence of this is that the limiting shape lies in a plane.

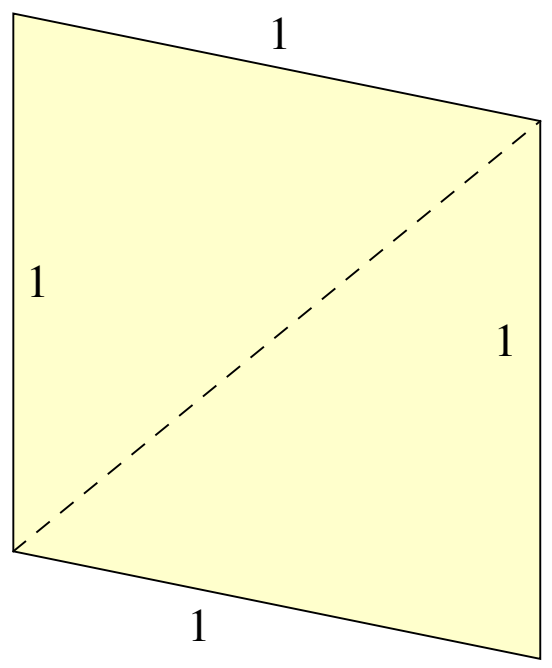

**Figure 13.** The normalized limiting rhombus; the shorter diagonal has length $\frac{2}{5}\sqrt{10}$.

Second, we can learn some interesting asymptotic behavior as follows. Let the notation $L \to^* \infty$ mean that $L$ takes on only values satisfying the condition of the proof: $\left|\overline{(L+1)\,\theta}\right| < \frac{2\pi}{\sqrt{5}\,(L+1)}$. Then $\lim_{L\to^*\infty} \frac{K-I}{\sin^2(\frac{\delta}{2})} = \frac{4}{3125} H_0$, which is the same as

$\lim_{L\to^*\infty} \frac{K-I}{\overline{\delta}^2} = \frac{1}{3125} H_0$. Also $\lim_{L\to^*\infty} \frac{H_0}{L} = 250 \begin{pmatrix} 3 & 3 & 3 & 3 \\ -1 & -1 & -1 & -1 \\ -1 & -1 & -1 & -1 \\ -1 & -1 & -1 & -1 \end{pmatrix}$. Therefore $\lim_{L\to\infty^*} \frac{K-I}{L\,\overline{\delta}^2} = \frac{2}{25} \begin{pmatrix} 3 & 3 & 3 & 3 \\ -1 & -1 & -1 & -1 \\ -1 & -1 & -1 & -1 \\ -1 & -1 & -1 & -1 \end{pmatrix}$. The norm of this last matrix is $\frac{8}{25}\sqrt{3}$, so we learn that, for nearly closed quadrahelices, $\|K-I\|_2$ is close to $\frac{8\sqrt{3}}{25} L\,\overline{\delta}^2$. Another way of looking at this is that the gap, as a fraction of the span of the whole chain, is infinitely often proportional to the square of $\overline{(L+1)\,\theta}$. The chart in Figure 14, based on all values of $L$ from 4 to 200000, illustrates the convergence to $8\sqrt{3}/25$.

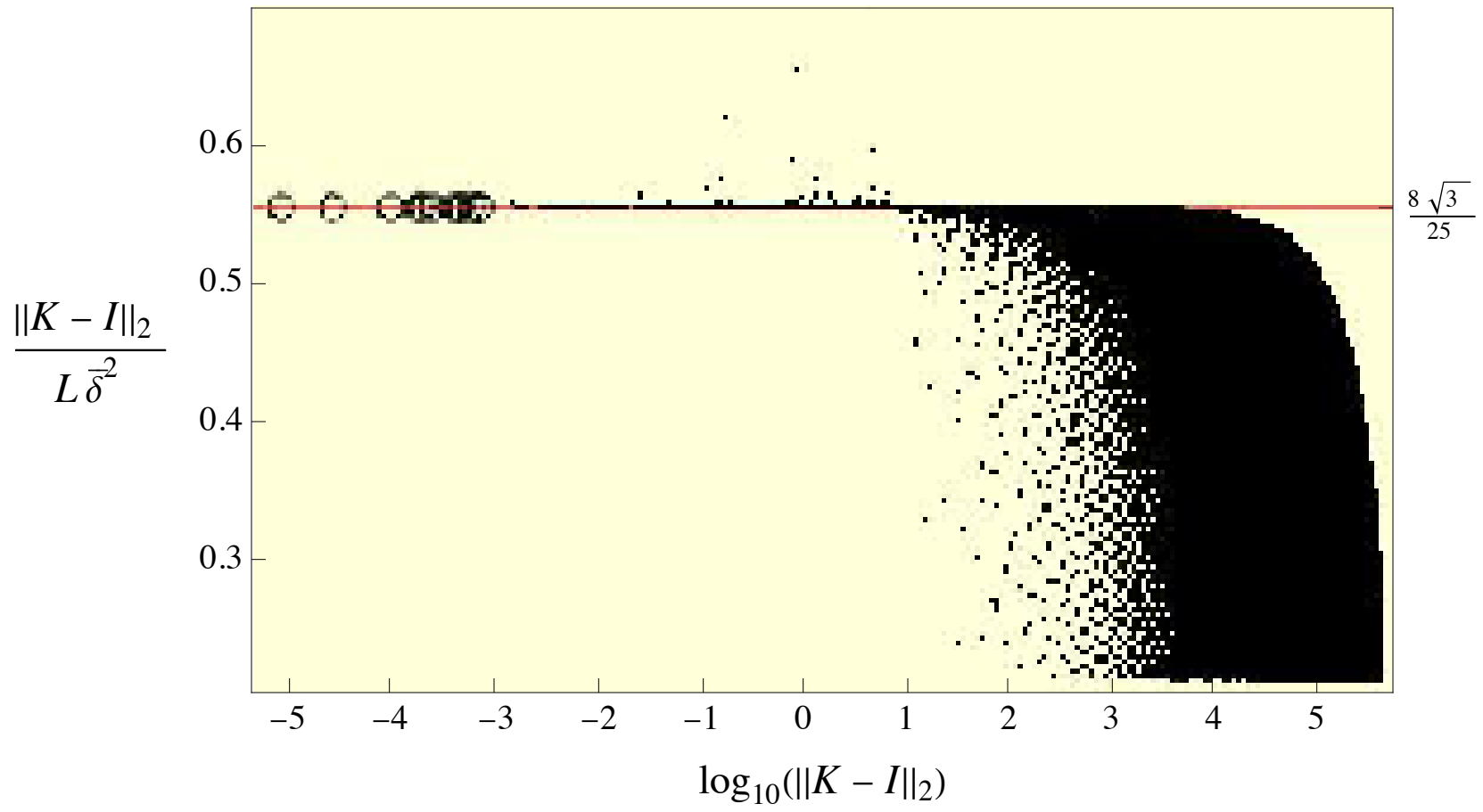

**Figure 14.** Let $K$ be the matrix defining the last tetrahedron in $QH_L$; $\overline{\delta}$ denotes $\overline{(L+1)\,\theta}$; The chart plots the ratio of $\|K-I\|_2$ to $L\,\overline{\delta}^2$ against $\log_{10}\|K-I\|_2$ for $L = 4, 5, \ldots, 200\,000$. Small norms are close to the product of $L\,\overline{\delta}^2$ and $\frac{8\sqrt{3}}{25}$.

The use of continued fractions quickly leads to nearly closed quadrahelices. Table 1 shows the convergents $\frac{k}{q}$ to $\frac{\theta}{2\pi}$ and the resulting small gaps in $QH_L$, where $L = q - 1$. This method easily gives an explicit chain with $10^{10^6}$ tetrahedra and correspondingly small gap. The table shows that $QH_{30\,170\,783\,468\,093\,193}$ has smaller gap than the best chain of [6] (shown in Fig. 16).

| $L = q - 1$ | $k = \text{round}\left[\frac{(L+1)\,\theta}{2\pi}\right]$ | $\overline{(L+1)\,\theta}$ | gap of $QH_L$ |
|---|---|---|---|
| 1 | 1 | −1.7 | 1.0 |
| 2 | 1 | 0.62 | 0.32 |
| 7 | 3 | −0.45 | 0.42 |
| 10 | 4 | 0.17 | 0.078 |
| 29 | 11 | −0.099 | 0.063 |
| 40 | 15 | 0.074 | 0.046 |
| 70 | 26 | −0.026 | 0.0095 |
| 182 | 67 | 0.022 | 0.018 |
| 253 | 93 | −0.0031 | 0.00050 |
| 1960 | 718 | 0.00048 | 0.000089 |
| 12 019 | 4401 | −0.00026 | 0.00016 |
| 13 980 | 5119 | 0.00022 | 0.00013 |
| 26 000 | 9520 | −0.000042 | $8.9\ 10^{-6}$ |
| 143 985 | 52 719 | 0.00001 | $2.9\ 10^{-6}$ |
| 601 944 | 220 396 | $-1.1\ 10^{-6}$ | $1.3\ 10^{-7}$ |
| 5 561 490 | 2 036 283 | $6.7\ 10^{-7}$ | $4.9\ 10^{-7}$ |
| 6 163 435 | 2 256 679 | $-3.8\ 10^{-7}$ | $1.8\ 10^{-7}$ |
| 30 170 783 468 093 193 | 30 170 783 468 093 193 | $1.3\times10^{-17}$ | $9.5\times10^{-19}$ |
| 52 126 933 878 205 ≪70≫ 925 856 156 787 510 | 19 085 743 247 326 ≪70≫ 845 985 782 583 721 | $1.6\ 10^{-100}$ | $2.7\ 10^{-102}$ |
| 1 114 768 425 712 ≪974≫ 60 678 210 425 128 | 1 114 768 425 712 ≪974≫ 60 678 210 425 128 | $1.5\times10^{-1000}$ | $4.9\times10^{-1001}$ |

**Table 1.** The left column shows the denominators of the convergents to $\frac{\theta}{2\pi}$, less 1; these yield almost closed chains $QH_L$. The second column has the numerators, the multiples of $2\pi$ that reduce $(L+1)\,\theta$ to near 0. The third column shows the reduced angle $\overline{(L+1)\,\theta}$ and the last gives the gap size, measured by Hausdorff distance.

# 5. Algebra and Geometry

An analysis of the matrix $K$ ties together several algebraic and geometric facts about a certain family of tetrahedral chains (which includes the quadrahelix and octahelix). Suppose $C$ is a chain of $n$ tetrahedra $T_i$ that is symmetric by reflection in a plane $\Pi$ through one of its faces; then $n$ is even and the reflection face is the middle face. Let $F$ denote reflection in $\Pi$; $F$ takes $T_1$ to $T_n$. Let $F'$ be the reflection in the face $\Pi'$ of $T_1$ that gives the invisible tetrahedron $T_0$. By Theorem 1, these two planes are not parallel. Let $\Omega$ be the isometry taking $T_0$ to $T_n$. Then $\Omega = F\,F'$ and therefore $\Omega$ is a rotation. The axis of $\Omega$ is the line $\Pi \cap \Pi'$ (Proof: for $p$ on this line, $\Omega(p) = F F'(p) = p$). This explains why our constructions have the useful rotation axis defined by certain planes (see Fig. 10). The chain is embedded if its first half is embedded and lives on one side of $\Pi$.

Although $T_i$ was defined to be the vertex set of a tetrahedron in a chain, we will also use it for the $3\times4$ matrix with these vertices as columns. With $K$ as the barycentric matrix of the chain (the product of the $M_i$), recall that the final tetrahedron $T_n$ is given by $T_n = T_0\,K$. This transformation of the invisible tetrahedron to the final tetrahedron can be viewed as a transformation of the $4\times4$ matrix $\mathcal{T}_0 = \begin{bmatrix} T_0 \\ 1\,1\,1\,1 \end{bmatrix}$ to $\mathcal{T}_n = \begin{bmatrix} T_n \\ 1\,1\,1\,1 \end{bmatrix}$. One way to effect this transformation is by $\mathcal{T}_0\,K = \mathcal{T}_n$. Another way is to consider the decomposition of $\Omega$, the rotation induced by the motion; $\Omega$ decomposes into $\tau \circ R$ where $R$ is a rotation that fixes the origin and $\tau$ is a translation by $\vec{t}$. Therefore, letting $\mathcal{R} = \begin{pmatrix} R & \vec{t} \\ 0\,0\,0 & 1 \end{pmatrix}$, we have $\mathcal{R}\,\mathcal{T}_0 = \mathcal{T}_n$. The two approaches combine to give $\mathcal{R}\,\mathcal{T}_0 = \mathcal{T}_0\,K$ and $\mathcal{R} = \mathcal{T}_0\,K\,\mathcal{T}_0^{-1}$. For $QH_L$, we derived a symbolic expression for $K$; for chains in general $K$ is easy to compute as a product. So the preceding relation gives both a symbolic expression and a simple method of computation for $R$ and $\vec{t}$.

The eigenvalues of $\mathcal{R}$ (and hence also $K$) are $z$, $\bar{z}$, 1, 1. The two left fixed points of $\mathcal{R}$ are $(0, 0, 0, 1)$ and $(w_1, w_2, w_3, 0)$ where, letting $\vec{w} = (w_1, w_2, w_3)$, $\vec{w}^T R = \vec{w}^T$ and $\vec{w} \cdot \vec{t} = 0$. The corresponding two fixed points of $K$ are $(0, 0, 0, 1)\,\mathcal{T}_0 = (1, 1, 1, 1)$ and $(w_1, w_2, w_3, 0)\,\mathcal{T}_0 = \vec{w}^T\,T_0$. One can get $\vec{w}$ as the normalized cross product of the largest (in norm) columns of $R - I$; so $\vec{w}$ is a unit vector.

We learn from the preceding that:

- $R$ and $\vec{t}$ are easily computed from $K$.

- $\vec{w}$ gives the direction of the axis of $R$ and $\Omega$ and is easily computed as an eigenvector of $R$.

- The rank of $K - I$ is 2 (it is not 1 because its eigenvalues are $(z, \bar{z}, 1, 1) - 1 = \left(z - 1, \overline{z - 1}, 0, 0\right)$. The left kernel of $K - I$ is generated by $(1, 1, 1, 1)$ and $\vec{w}^T\, T_0$.

To fully characterize the isometry $\Omega$ we need a point on its axis; such a point $\vec{u}$ satisfies $\vec{u} - R\,\vec{u} = \vec{t}$. But $R - I$ has rank 2 so is singular (as is true of any $R \in SO_3(\mathbb{R})$), so we cannot get $\vec{u}$ by inversion. But fairly standard linear algebra gives $\vec{u}$ as follows. Let $\vec{t}^*$ be $\vec{t}$ normalized to have length 1. Then $\vec{u} = Q\left(I_2 - Q^T R\, Q\right)^{-1} Q^T\, \vec{t}$, where $Q$ is the $3\times 2$ matrix whose first column is $\vec{t}^*$ and second column is $\vec{w} \times \vec{t}^*$. So now $\Omega$ is completely characterized: $\Omega$ is a rotation around the axis in direction $\vec{w}$ through the point $\vec{u}$ and by an angle that is $\arg(v)$ where $v$ is an eigenvalue of $R$; thus the line $\Pi \cap \Pi'$ is given by $\vec{u} + \alpha\, \vec{w}$.

These general ideas can be applied to the quadrahelix and lead to a different proof of Theorem 5, one that has much in common with the geometric proof mentioned prior to that theorem.

- By Corollary 3, $\|R - I\| \le \bar{\delta}^2$, where $\bar{\delta}$ is the mod-$2\pi$ reduction of $(L + 1)\,\theta$.

- The proof of the acute angle lemma showed that $\rho_0 \ge \cos^{-1}\left(\frac{2\sqrt{6}}{5}\right)$; applying the Law of Cosines to an isosceles triangle with apex angle $2\rho_0$ and opposite side bounded by $(2\,L + 1)\,h$ yields

$$\left\|\vec{u}\right\| \le \frac{h\,(2\,L+1)}{\sqrt{2}\,\sqrt{1 - \cos\left(2\cos^{-1}\left(\frac{2\sqrt{6}}{5}\right)\right)}} = \frac{h\,(2\,L+1)}{\sqrt{2}\,\sqrt{1 - \frac{23}{25}}} = \frac{5}{2}\,h\,(2\,L + 1) \le 0.79 + 1.59\,L$$

If we assume $L \ge 17$, then $1 + \left\|\vec{u}\right\| \le 1.7\,L$.

- $\|\mathcal{T}_0\| = \frac{1}{5\sqrt{2}}\sqrt{117 + \sqrt{8689}} \le 2.06$

- $\text{gap} \le \|K - I\| \le \frac{1}{\sigma_{\min}(T_0)}\,\|T_0\,(K - I)\| = \sqrt{2}\,\|T_n - T_0\|$, where $\sigma_{\min}$ denotes the smallest singular value.

- Let $\mathcal{R}_{3,4} = \left[(R - I)\quad \vec{t}\right]$, where $\left[\cdot\quad \vec{t}\right]$ indicates the adjunction of the column $\vec{t}$; then $T_n - T_0 = \mathcal{R}_{3,4}\,\mathcal{T}_0$,

These facts and an application of the induced 2-norm to the last equality give the following, assuming $L \ge 17$.

$$\begin{aligned}\frac{1}{\sqrt{2}}\,\text{gap}(\text{OH}_L) &\le \|T_n - T_0\| = \|\mathcal{R}_{3,4}\ \mathcal{T}_0\| \le \left\|\left[(R - I)\quad -(R - I)\,\vec{u}\right]\right\|\,\|\mathcal{T}_0\| \\ &\le \|R - I\|\left(1 + \left\|\vec{u}\right\|\right)2.06 \le \bar{\delta}^2\ (1.7\,L)\,2.06\ \le 3.51\,L\,\bar{\delta}^2\end{aligned}$$

So this bounds the gap by $5\,L\,\bar{\delta}^2$. Because $\bar{\delta}$ is less than $1/L$ infinitely often (see proof of Thm. 5) this provides an alternate proof that the gap of $\text{QH}_L$ can be made arbitrarily small.

# 6. An Octagonal Pattern

We first resolved the vanishing-gap conjecture with an 8-sided shape that is more complicated than the quadrahelix. The *octahelix* $\text{OH}_L$ arises from the 8-part sequence $S_{L+1}\ \overline{S}_L\ p(S_{L+1})\ p(\overline{S}_L)\ S_{L+1}\ \overline{S}_L\ p(S_{L+1})\ p(\overline{S}_L)$, where $S_m$ is the $m$-term string 12341234… and $p$ is the permutation $\{2, 3, 4, 1\}$, a 4-cycle. Thus $\text{OH}_4$ is 12341 4321 23412 1432 12341 4321 23412 1432. The useful symmetries of the shape are clarified if we shift the string left $L$ units; that is, use $j\ \overline{S}_L\ p\left(S_{L+1}\ \overline{S}_L\right) S_{L+1}\ \overline{S}_L\ p\left(S_{L+1}\ \overline{S}_L\right) S_L$ for $\text{OH}_L$, where $j \equiv L + 1 \pmod 4$. The cases $L = 4, 5, 6$, and 36 are shown in Figure 15.

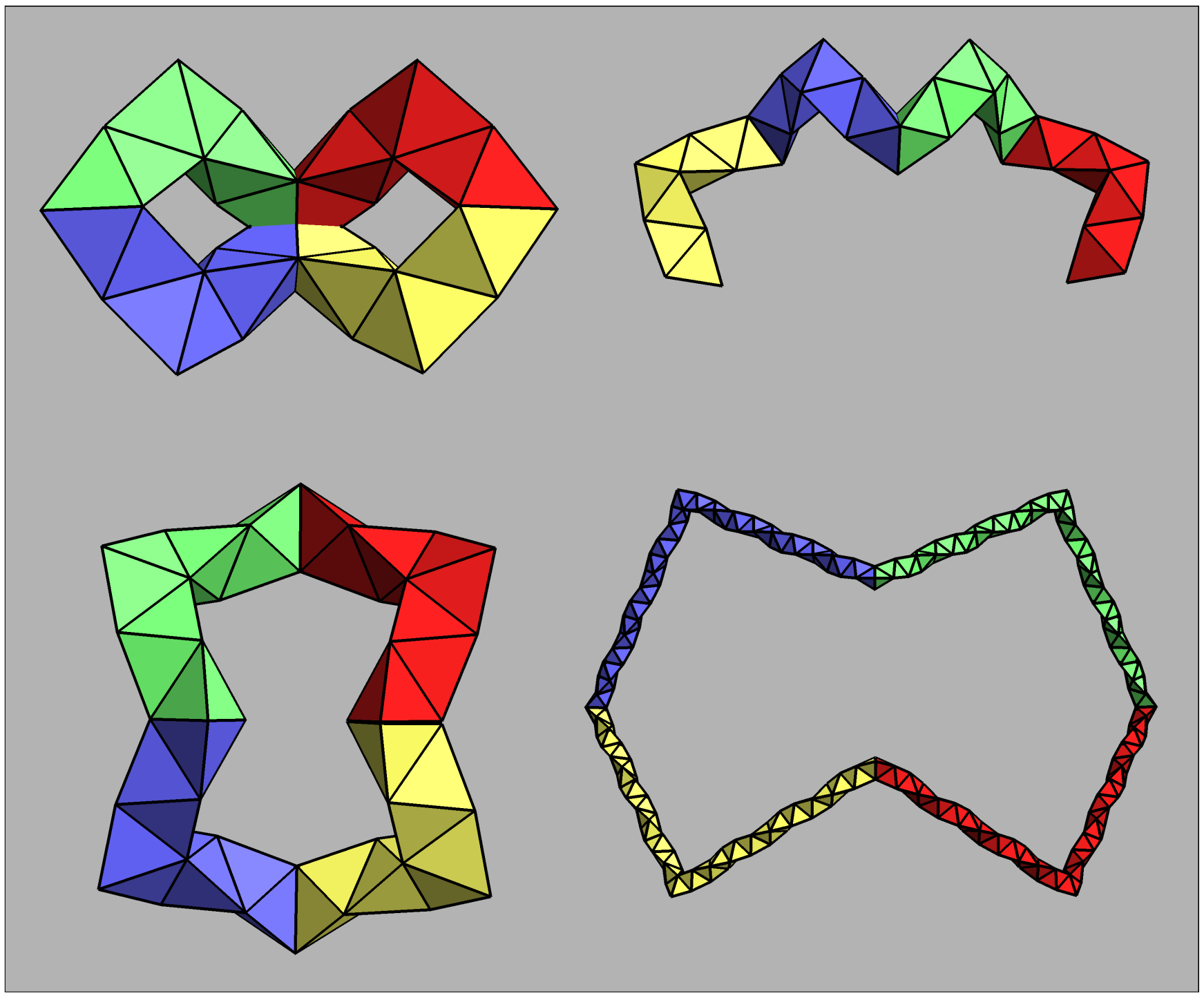

**Figure 15.** The octahelix for $L = 4, 5, 6$, and 36; $OH_4$ is not embedded. The gaps for $OH_6$ and $OH_{36}$ are both about 0.02.

The analysis is more complicated than for the quadrahelix, but the general approach is similar: there are five important planes (really nine, since the planes bisecting each colored segment are useful too), an acute angle lemma, a five-plane lemma, an embedding theorem (which fails for $L = 1$ or 4), two magic angles, and a number-theoretic relation that leads to small gaps. There are two magic angles, $\gamma^{\pm} = \cos^{-1}\left(\frac{1}{12}\left(-3 \pm 5\sqrt{3}\right)\right)$: when $\overline{L\theta}$ is near either one, the gap is small. If $\delta = L\theta - \gamma^{\pm}$ and assuming $|\delta| < 1$, one can show that the gap is bounded by $3L\overline{\delta}^2$. For example, if $L = 686$, then $\overline{686\,\theta} = \gamma^{+} + 0.00035\ldots$, the bound is 0.0003, and the actual gap in $OH_{686}$ is about 0.000016. A pleasant property of the octahelix is that all seven angles equal $\sec^{-1}(5)$.

The octahelix construction is quite similar to the quadrahelix, except that for QH there was only one magic angle, $-\theta$. When $\overline{L\theta}$ is sufficiently close to $-\theta$, the quadrahelical gap is small. Having $\theta$ be both the magic angle and the multiplying angle is wonderful because $L\theta + 2\pi k \sim -\theta$ becomes $(L+1)\theta + 2\pi k \sim 0$ and the continued fraction of $\frac{\theta}{2\pi}$ gives nearly closed quadrahelices. For the octahelix, we need the more complicated inhomogeneous relation $L\theta + 2\pi k \sim \gamma^{\pm}$. Kronecker's theorem and the density in the unit circle of the set $\{\overline{L\theta}\}$ can be used as mentioned in §4 and yields a proof that for a subsequence of $L$-values, $OH_L$ is embedded and has vanishingly small gap. But to get explicit $L$-values requires more sophistication than the continued fraction method. Methods of lattice reduction can be used to find explicit values, though it is not guaranteed to work in general. We present the details because the algorithm is interesting and useful. It works in our specific case, quickly providing values of $L$ for which the gap in $OH_L$ is $10^{-100}$ or smaller. We use $\lfloor x \rceil$ for the integer closest to $x$.

## Lattice Reduction method for the Diophantine relation $x\alpha + y\beta + \gamma \sim 0$

We have $\theta = \cos^{-1}\left(-\frac{2}{3}\right)$ and $\gamma^{\pm} = \cos^{-1}\left(\frac{1}{12}\left(-3 \pm 5\sqrt{3}\right)\right)$ and seek $L$ so that $\overline{L\theta}$ is close to $\gamma^{+}$ (or $\gamma^{-}$), with error about $\frac{1}{L}$. This is the same as $L\theta + 2\pi k$ being close to $\gamma^{+}$, which is a special case of this more general problem: Given $\alpha, \beta, \gamma \in \mathbb{R}$ with irrational $\frac{\alpha}{\beta}$, find $x, y, z \in \mathbb{Z}$ so that $x\alpha + y\beta - \gamma$ is close to 0 and $|x|$ is bounded by roughly the product of $|\beta|$ and the reciprocal of the error. Such integers do exist: Kronecker's Theorem [10, Thm. 440] gives, for any $X$, integers $x, y$ so that $x \geq X$ and $\left|x\frac{\alpha}{\beta} + y - \frac{\gamma}{\beta}\right| < \frac{3}{x}$, or $|x\alpha + y\beta - \gamma| < 3\left|\frac{\beta}{x}\right|$. In our problem $\alpha = \theta$, $\beta = 2\pi$, and $\gamma = \gamma^{+}$.

It is well-known that lattice reduction (see [3]) can provide a heuristic method for finding the approximating integers guaranteed by Kronecker's Theorem; a detailed discussion is in [9]. The idea is based on Babai's nearest plane algorithm. Here we take the point of view: Given $X$, can we find the values $x$, $y$ promised by Kronecker?

**Heuristic Algorithm for Diophantine Approximation**

*Input*. Reals $\alpha, \beta, \gamma$, and positive $X$, where $\alpha/\beta$ is irrational.

*Output*. Integers $x$ and $y$ that *might* satisfy $|x| \geq X$ and $|x\,\alpha + y\,\beta - \gamma| < 3\left|\frac{\beta}{x}\right|$

1. Initialize: Let $\epsilon = \left(3\left|\frac{\beta}{x}\right|\right)^2$; set a multiplicative factor $s = \epsilon^{-2}\max(1, |\alpha|, |\beta|)$; get integer approximations $a = \lfloor s\,\alpha\rceil$, $b = \lfloor s\,\beta\rceil$, $c = \lfloor s\,\gamma\rceil$; set $t = (-c, 0, 0)$.

2. Apply Lenstra–Lenstra–Lovász lattice reduction to the lattice generated by $(a, 1, 0)$ and $(b, 0, 1)$ to get $\{B_1, B_2\}$.

3. Orthogonalize to get $B^* = \left\{B_1, B_2 - \frac{B_2 \cdot B_1}{B_1 \cdot B_1} B_1\right\}$.

4. Subtract approximation to component in $B_2$ direction: $t = t - \left\lfloor \frac{t \cdot B_2^*}{B_2^* \cdot B_2^*} \right\rceil B_2$.

5. Subtract approximation to component in $B_1$ direction: $t = t - \left\lfloor \frac{t \cdot B_1^*}{B_1^* \cdot B_1^*} \right\rceil B_1$.

6. The desired coefficients are $x = t_2$ and $y = t_3$.

| $x$ | $y$ | $\log_{10}$ error |
|---:|---:|---:|
| 4 | −1 | −1.80 |
| 686 | −251 | −3.46 |
| 1274 | −466 | −3.88 |
| 64 708 | −23 692 | −5.48 |
| 666 653 | −244 088 | −5.65 |
| 1 870 543 | −684 880 | −6.86 |
| 111 021 125 | −40 649 248 | −7.79 |
| 233 817 317 | −85 609 817 | −9.23 |
| 3 113 400 370 | −1 139 939 675 | −9.71 |
| 434 337 601 428 | −159 028 266 709 | −10.4 |

**Table 2.** The output of the Diophantine approximation method: $x\,\theta$, reduced mod $2\pi$, is very close to one of $\gamma^{\pm}$. The errors shown are always less than $6\,\pi/x$, consistent with Kronecker’s theorem.

For our octahelix case involving $\theta$, $2\pi$, and $\gamma^+$, the algorithm might return a negative integer $x$. When that happens it is easy to see that $-x-1$ yields small error when the conjugate value $\gamma^-$ is used instead of $\gamma^+$. We have made that change in Table 2 for the cases where $x$ was negative. Using this method it takes only an instant to find that for the 101-digit integer $L =$

10 021 748 859 140 317 070 670 606 276 035 026 550 890 854 358 977 791 419 357 055 984 002 968 975 063 788 933 115 - 779 224 144 920 189

the difference between $\overline{L\,\theta}$ and $\gamma^-$ is less than $10^{-100}$, with correspondingly small gap in $\mathrm{OH}_L$. For the general problem, a few thousand random trials for irrational values $\alpha, \beta, \gamma$ show that the heuristic method seems to work quite well.

## 7. Open Questions

A natural question is whether our quadrahelix can be improved. There are several aspects to consider.

**Question 1.** Is there a triangular pattern that leads to embedded chains with vanishing gaps?

Some investigations indicate that the quadrahelix idea will not work for a pattern based on an equilateral triangle or a square. But there might be something involving other types of triangles, or a rectangle, or perhaps a completely sporadic sequence of triangles that works.

The quadrahelix requires roughly $6/\epsilon$ tetrahedra to achieve a gap of $\epsilon$; it takes about $10^{17}$ tetrahedra to get a quadrahelical gap near $10^{-17}$. In [6] we found an embedded chain of 540 tetrahedra with a gap of $7\cdot 10^{-18}$ (see Fig. 16). The reflection sequence for this chain is $\left((b\,4\,\overline{b})\;(b\,4\,\overline{b})_{\pi}\;(b\,4\,\overline{b})\,(b\,4\,\overline{b})_{\pi^3}\right)^3$ where $b = 1\,234\,123\,434\,132\,341\,213\,412$ (length 22) and $\pi$ is the 4-cycle $\{2, 3, 4, 1\}$; total length is $(2\cdot 22+1)\cdot 4\cdot 3 = 540$. Experiments [6, Fig. 9] suggest that a gap of $\epsilon$ can be achieved with $10\log(1/\epsilon)$ tetrahedra, a formula that predicts about 400 tetrahedra for $\epsilon = 10^{-17}$.

**Question 2.** Is there a pattern for embedded chains that achieves gap $\epsilon$ using $O(\log(1/\epsilon))$ tetrahedra?

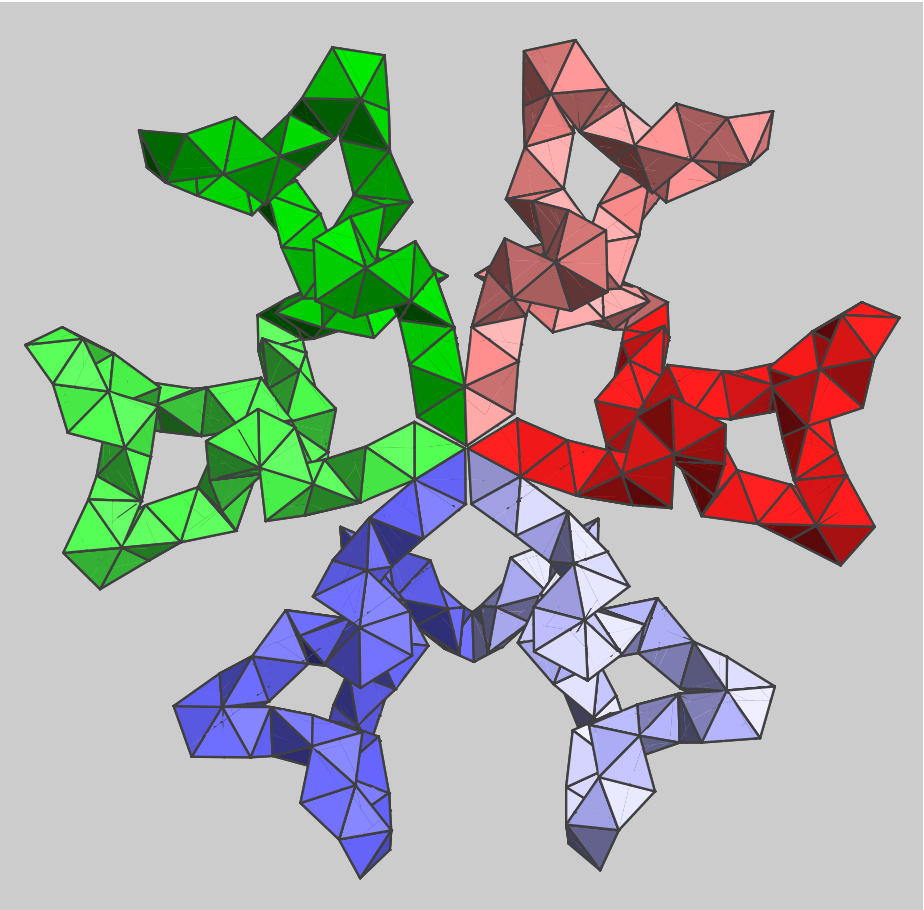

**Figure 16.** An embedded chain of 540 tetrahedra that has a gap smaller than $10^{-17}$.

A main point about our symmetric patterns is that the gap depends on two things: how close the rotation angle induced by the chain is to $2\pi$, and how large $L$ is. More generally, it would seem that the role of $L$ might be replaced by the overall diameter of the chain, so that one line of attack on Question 2 is to look for chains that use a large number of tetrahedra but have them arranged much more compactly than the quadrahelix or octahelix.

# Appendix: Formulas and *Mathematica* Code

## Basic Utilities

```
Clear[a, b, c, x1, y1, z1, x2, y2, z2, x3, y3, z3, x4, y4, z4];
bary[{{x1_, y1_, z1_}, {x2_, y2_, z2_}, {x3_, y3_, z3_}, {x4_, y4_, z4_}}, {a_, b_, c_}] :=
  Evaluate[With[{pts4 = {{x1, y1, z1}, {x2, y2, z2}, {x3, y3, z3}, {x4, y4, z4}}, pt = {a, b, c}, var = x /@ Range[4]},
    Simplify[var /. Solve[{Total[var] == 1, And @@ Thread[var.pts4 == pt]}, var][[1]]]]];
fromBary[bary_, pts4_] := bary.pts4;
θ = ArcCos[-2/3]; r = 3 Sqrt[3]/10; h = 1/Sqrt[10]; id = IdentityMatrix[4];
V[i_] := {r Cos[i θ], r Sin[i θ], i h};
base1 = V /@ {-1, 0, 1, 2};
```

## The Barycentric Formula for the Tetrahelix

```
CC = Simplify[bary[base1, V[Q]]] /. {Cos[Q θ] → c, Sin[Q θ] → s};
{coeffQ, coeffc, coeffs} = Coefficient[CC, #] & /@ {Q, c, s}
con = Expand[CC - (coeffQ Q + coeffc + coeffs s)]
```

$$\left\{\left\{-\frac{3}{10}, -\frac{1}{10}, \frac{1}{10}, \frac{3}{10}\right\}, \left\{-\frac{3}{10}, \frac{3}{5}, -\frac{3}{10}, 0\right\}, \left\{-\frac{3}{5\sqrt{5}}, \frac{3}{10\sqrt{5}}, \frac{6}{5\sqrt{5}}, -\frac{9}{10\sqrt{5}}\right\}\right\}$$

$$\left\{\frac{3}{5}-\frac{3c}{10}, -\frac{1}{5}+\frac{3c}{5}, \frac{3}{5}-\frac{3c}{10}, 0\right\}$$

## The Dot Product in the Acute Angle Lemma

```
A = TrigExpand[(3 V[3] - Total[{V[0], V[1], V[2]}])/Sqrt[6]]
3 Sqrt[5] A
```

$$\left\{\frac{\sqrt{2}}{3}, \frac{2\sqrt{\frac{2}{5}}}{3}, \sqrt{\frac{3}{5}}\right\}$$

$$\left\{\sqrt{10}, 2\sqrt{2}, 3\sqrt{3}\right\}$$

```
B = Together[TrigExpand[V[L + 3] - V[L + 1] /. f_[z_] /; f === Sin || f === Cos :> f[Expand[z]]]];
TrigExpand[ Sqrt[15] B /. θ L → δ - θ]
```

$$\left\{-\sqrt{5}\,\mathrm{Cos}[\delta]+2\,\mathrm{Sin}[\delta], -2\,\mathrm{Cos}[\delta]-\sqrt{5}\,\mathrm{Sin}[\delta], \sqrt{6}\right\}$$

```
Factor[TrigExpand[A.B /. L θ → δ - θ]]
TrigFactor[TrigExpand[A.B /. L θ → δ - θ]]
```

$$-\frac{1}{5}\sqrt{6}\,(-1+\mathrm{Cos}[\delta])$$

$$\frac{2}{5}\sqrt{6}\,\mathrm{Sin}\left[\frac{\delta}{2}\right]^2$$

## The Determinant for the Embedding Theorem

```
mat = TrigExpand[(V[Q] - # &) /@ {V[1], V[2], 1/2 (V[0] + V[3])}];
MatrixForm[mat /. θ → "θ"]
X = Expand[20 √10 Det[mat]];
X /. θ → "θ"
```

$$\begin{pmatrix} \frac{\sqrt{3}}{5} + \frac{3}{10}\sqrt{3}\,\mathrm{Cos}[\theta Q] & -\frac{\sqrt{\frac{3}{5}}}{2} + \frac{3}{10}\sqrt{3}\,\mathrm{Sin}[\theta Q] & -\frac{1}{\sqrt{10}} + \frac{Q}{\sqrt{10}} \\ \frac{1}{10\sqrt{3}} + \frac{3}{10}\sqrt{3}\,\mathrm{Cos}[\theta Q] & \frac{2}{\sqrt{15}} + \frac{3}{10}\sqrt{3}\,\mathrm{Sin}[\theta Q] & -\sqrt{\frac{2}{5}} + \frac{Q}{\sqrt{10}} \\ -\frac{49}{60\sqrt{3}} + \frac{3}{10}\sqrt{3}\,\mathrm{Cos}[\theta Q] & -\frac{7}{12\sqrt{15}} + \frac{3}{10}\sqrt{3}\,\mathrm{Sin}[\theta Q] & -\frac{3}{2\sqrt{10}} + \frac{Q}{\sqrt{10}} \end{pmatrix}$$

$$-9\sqrt{5} + 6\sqrt{5}\,Q - \sqrt{5}\,\mathrm{Cos}[\theta Q] + 7\,\mathrm{Sin}[\theta Q]$$

```
N[X /. z_ (Cos | Sin)[_] :→ -Abs[z]]
```

-29.3607 + 13.4164 Q

## The Symbolic Matrix *K* of the Quadrahelix

First get the barycentric form of a general point $V_Q$ on the tetrahelix. It is not too complicated.

```
baryGen = FullSimplify[bary[base1, V[Q]]];
Column[Expand[baryGen /. {Cos[Q θ] → c, Sin[Q θ] → s}]]
```

$$\frac{3}{10} - \frac{3c}{10} - \frac{3Q}{10} - \frac{3s}{5\sqrt{5}}$$
$$\frac{2}{5} + \frac{3c}{5} - \frac{Q}{10} + \frac{3s}{10\sqrt{5}}$$
$$\frac{3}{10} - \frac{3c}{10} + \frac{Q}{10} + \frac{6s}{5\sqrt{5}}$$
$$\frac{3Q}{10} - \frac{9s}{10\sqrt{5}}$$

```
base2 = Table[V[i], {i, L + {4, 1, 2, 3}}];
base3 = Table[fromBary[baryGen, base2], {Q, L + {3, 2, 1, 4}}];
base4 = Table[fromBary[baryGen, base3], {Q, L + {3, 0, 1, 2}}];
finalTetBaryForm =
  Table[bary[base1, Simplify[TrigExpand[ExpandAll[fromBary[baryGen, base4]]]]], {Q, L + {0, 3, 2, 1}}];
matrixK = Transpose[finalTetBaryForm];

(*raw=Collect[Together[
    ExpandAll[3125 TrigExpand[matrixK-id/.L θ→δ-θ]/.{Cos[δ]→c,Sin[δ]→s}/.{c^2→1-s^2,c^3→c (1-s^2),c^4→(1-s^2)^2}/.
      {c→1-2 σ,s^2→4 σ-4 σ^2,s^3→s (4 σ-4 σ^2),s^4→(4 σ-4 σ^2)^2}]],σ];*)

raw =
  ExpandAll[3125 TrigExpand[matrixK - id /. L θ → δ - θ] /. {Cos[δ] → c, Sin[δ] → s} /. {c^z_ :→ (1 - s^2)^Floor[z/2] c^Mod[z,2]} /.
    {c → 1 - 2 σ, s^z_ :→ (4 σ - 4 σ^2)^Floor[z/2] s^Mod[z,2]}];
```

Next is therefore an exact symbolic form of the barycentric matrix $K$ for $\mathrm{QH}_L$, where $s = \sin\delta$, $\sigma = \sin^2\left(\frac{\delta}{2}\right)$ and $\delta = (L+1)\,\theta$.

```
1/3125 raw + id // Simplify // Column
```

$$\Big\{\frac{3125+300\left(10+10\,L-\sqrt{5}\,s\right)\sigma-24\,000\,\sigma^2-576\left(10+10\,L-\sqrt{5}\,s\right)\sigma^3+23\,040\,\sigma^4}{3125},$$

$$\frac{12\,\sigma\left(-25+48\,\sigma^2\right)\left(-2\sqrt{5}\,L+6\,s+\sqrt{5}\,(-7+6\,\sigma)\right)}{625\sqrt{5}},\ -\frac{12\left(2\sqrt{5}+2\sqrt{5}\,L+9\,s\right)\sigma\left(-25+48\,\sigma^2\right)}{625\sqrt{5}},\ -\frac{12\,\sigma\left(-25+48\,\sigma^2\right)\left(2\sqrt{5}\,L-6\,s+3\sqrt{5}\,(-1+2\,\sigma)\right)}{625\sqrt{5}}\Big\}$$

$$\Big\{-\frac{4\,\sigma\left(3\,s\left(-175-120\,\sigma+384\,\sigma^2\right)+L\left(50\sqrt{5}+240\left(\sqrt{5}-s\right)\sigma-384\sqrt{5}\,\sigma^2\right)+2\sqrt{5}\left(275-380\,\sigma-624\,\sigma^2+720\,\sigma^3\right)\right)}{625\sqrt{5}},\ \frac{1}{625\sqrt{5}}$$

$$\left(625\sqrt{5}-100\left(7\sqrt{5}+2\sqrt{5}\,L-6\,s\right)\sigma-120\left(23\sqrt{5}+8\,L\left(\sqrt{5}-s\right)-52\,s\right)\sigma^2+192\left(31\sqrt{5}+8\sqrt{5}\,L-39\,s\right)\sigma^3-2304\sqrt{5}\,\sigma^4\right),$$

$$-\frac{4\,\sigma\left(s\left(225+840\,\sigma-1728\,\sigma^2\right)+L\left(50\sqrt{5}+240\left(\sqrt{5}-s\right)\sigma-384\sqrt{5}\,\sigma^2\right)+2\sqrt{5}\left(25+120\,\sigma-624\,\sigma^2+432\,\sigma^3\right)\right)}{625\sqrt{5}},$$

$$\frac{4\,\sigma\left(L\left(-50\sqrt{5}-240\left(\sqrt{5}-s\right)\sigma+384\sqrt{5}\,\sigma^2\right)+3\left(s\left(50+120\,\sigma-144\,\sigma^2\right)+\sqrt{5}\left(25+70\,\sigma-624\,\sigma^2+576\,\sigma^3\right)\right)\right)}{625\sqrt{5}}\Big\}$$

$$\Big\{-\frac{4\,\sigma\left(s\left(975+480\,\sigma-1872\,\sigma^2\right)+2\sqrt{5}\left(25-200\,\sigma-144\,\sigma^2+288\,\sigma^3\right)+L\left(480\,s\,\sigma+2\sqrt{5}\left(25-48\,\sigma^2\right)\right)\right)}{625\sqrt{5}},$$

$$-\frac{4\,\sigma\left(-6\,s\left(25-280\,\sigma+192\,\sigma^2\right)+\sqrt{5}\left(175-150\,\sigma-1488\,\sigma^2+1440\,\sigma^3\right)+L\left(480\,s\,\sigma+2\sqrt{5}\left(25-48\,\sigma^2\right)\right)\right)}{625\sqrt{5}},$$

$$\frac{625\sqrt{5}-100\left(2\sqrt{5}+2\sqrt{5}\,L+9\,s\right)\sigma-1920\,(1+L)\,s\,\sigma^2+192\left(-34\sqrt{5}+2\sqrt{5}\,L+9\,s\right)\sigma^3+6912\sqrt{5}\,\sigma^4}{625\sqrt{5}},$$

$$-\frac{4\,\sigma\left(L\left(480\,s\,\sigma+2\sqrt{5}\left(25-48\,\sigma^2\right)\right)+3\left(s\left(-50-240\,\sigma+576\,\sigma^2\right)+\sqrt{5}\left(-25+50\,\sigma-336\,\sigma^2+288\,\sigma^3\right)\right)\right)}{625\sqrt{5}}\Big\}$$

$$\Big\{\frac{4\,\sigma\left(s\left(525+120\,\sigma-864\,\sigma^2\right)+L\left(-50\sqrt{5}+240\left(\sqrt{5}+s\right)\sigma-192\sqrt{5}\,\sigma^2\right)+2\sqrt{5}\left(225+20\,\sigma-624\,\sigma^2+432\,\sigma^3\right)\right)}{625\sqrt{5}},$$

$$\frac{4\,\sigma\left(-6\,s\left(-25-20\,\sigma+24\,\sigma^2\right)+L\left(-50\sqrt{5}+240\left(\sqrt{5}+s\right)\sigma-192\sqrt{5}\,\sigma^2\right)+\sqrt{5}\left(-175+990\,\sigma-1968\,\sigma^2+1152\,\sigma^3\right)\right)}{625\sqrt{5}},$$

$$-\frac{4\,\sigma\left(3\,s\left(75-440\,\sigma+288\,\sigma^2\right)+2\,L\left(25\sqrt{5}-120\left(\sqrt{5}+s\right)\sigma+96\sqrt{5}\,\sigma^2\right)+2\sqrt{5}\left(25-120\,\sigma-336\,\sigma^2+432\,\sigma^3\right)\right)}{625\sqrt{5}},$$

$$\frac{625\sqrt{5}-100\left(2\sqrt{5}\,L-3\left(\sqrt{5}+2\,s\right)\right)\sigma+120\left(-17\sqrt{5}-36\,s+8\,L\left(\sqrt{5}+s\right)\right)\sigma^2-192\left(4\sqrt{5}\,L-9\left(\sqrt{5}+3\,s\right)\right)\sigma^3}{625\sqrt{5}}\Big\}$$

# A Compact Form of *K*

```
raw1 = Expand[raw/(4 σ)];
cσ = CoefficientList[raw1, σ];

Do[cL[k] = Table[Coefficient[cσ⟦i, j, k + 1⟧, L], {i, 4}, {j, 4}];
 cs[k] = Table[Coefficient[cσ⟦i, j, k + 1⟧, s], {i, 4}, {j, 4}];
 cCon[k] = Simplify[Table[cσ⟦i, j, k + 1⟧ - cL[k]⟦i, j⟧ L - cs[k]⟦i, j⟧ s, {i, 4}, {j, 4}]], {k, 0, 2, 2}]

cL[3] = Table[Coefficient[If[Length[cσ⟦i, j⟧] == 4, cσ⟦i, j, 4⟧, 0], L], {i, 4}, {j, 4}];
cs[3] = Table[Coefficient[If[Length[cσ⟦i, j⟧] == 4, cσ⟦i, j, 4⟧, 0], s], {i, 4}, {j, 4}];
cCon[3] = Simplify[Table[If[Length[cσ⟦i, j⟧] == 4, cσ⟦i, j, 4⟧, 0] - cL[3]⟦i, j⟧ L - cs[3]⟦i, j⟧ s, {i, 4}, {j, 4}]];

cσTemp = cσ /. L s → crossTerm;
cL[1] = Table[Coefficient[cσTemp⟦i, j, 2⟧, L], {i, 4}, {j, 4}];
cs[1] = Table[Coefficient[cσTemp⟦i, j, 2⟧, s], {i, 4}, {j, 4}];
cCross[1] = Table[Coefficient[cσTemp⟦i, j, 2⟧, crossTerm], {i, 4}, {j, 4}];
cCon[1] =
  Simplify[Table[cσTemp⟦i, j, 2⟧ - (cL[1]⟦i, j⟧ L + cs[1]⟦i, j⟧ s + cCross[1]⟦i, j⟧ crossTerm), {i, 4}, {j, 4}]];
```

```
Row[{125 MatrixForm[cCon[0]/125], " + ", 250 MatrixForm[cL[0]/250], "L", " + ", 75 √5 MatrixForm[cs[0]/(75 √5)], "sin(δ)"}]
Row[{50 MatrixForm[cCon[1]/50], " + ", 1200 MatrixForm[cL[1]/1200], "L", " + ",
  6 √5 MatrixForm[cs[1]/(6 √5)], "sin(δ)", " + ", 240 √5 MatrixForm[cCross[1]/(240 √5)], "L sin(δ)"}]
Row[{5 48 MatrixForm[cCon[2]/(5 48)], " + ", 480 MatrixForm[cL[2]/480], "L", " + ", 48 3 √5 MatrixForm[cs[2]/(48 3 √5)], "sin(δ)"}]
Row[{1440 MatrixForm[cCon[3]/1440], " + ", 0 MatrixForm[cL[3]/10], " L", " + ", 0 MatrixForm[cs[3]/(3 √5)], " sin(δ)"}]
```

$$125 \begin{pmatrix} 6 & 21 & 6 & -9 \\ -22 & -7 & -2 & 3 \\ -2 & -7 & -2 & 3 \\ 18 & -7 & -2 & 3 \end{pmatrix} + 250 \begin{pmatrix} 3 & 3 & 3 & 3 \\ -1 & -1 & -1 & -1 \\ -1 & -1 & -1 & -1 \\ -1 & -1 & -1 & -1 \end{pmatrix} L + 75\sqrt{5} \begin{pmatrix} -1 & -6 & 9 & -6 \\ 7 & 2 & -3 & 2 \\ -13 & 2 & -3 & 2 \\ 7 & 2 & -3 & 2 \end{pmatrix} \sin(\delta)$$

$$50 \begin{pmatrix} -120 & -45 & 0 & 45 \\ 76 & -69 & -24 & 21 \\ 40 & 15 & 0 & -15 \\ 4 & 99 & 24 & -51 \end{pmatrix} + 1200 \begin{pmatrix} 0 & 0 & 0 & 0 \\ -1 & -1 & -1 & -1 \\ 0 & 0 & 0 & 0 \\ 1 & 1 & 1 & 1 \end{pmatrix} L + 6\sqrt{5} \begin{pmatrix} 0 & 0 & 0 & 0 \\ 60 & 260 & -140 & 60 \\ -80 & -280 & -80 & 120 \\ 20 & 20 & 220 & -180 \end{pmatrix} \sin(\delta) + 240\sqrt{5} \begin{pmatrix} 0 & 0 & 0 & 0 \\ 1 & 1 & 1 & 1 \\ -2 & -2 & -2 & -2 \\ 1 & 1 & 1 & 1 \end{pmatrix} L \sin(\delta)$$

$$240 \begin{pmatrix} -6 & -21 & -6 & 9 \\ 26 & 31 & 26 & -39 \\ 6 & 31 & -34 & 21 \\ -26 & -41 & 14 & 9 \end{pmatrix} + 480 \begin{pmatrix} -3 & -3 & -3 & -3 \\ 4 & 4 & 4 & 4 \\ 1 & 1 & 1 & 1 \\ -2 & -2 & -2 & -2 \end{pmatrix} L + 144\sqrt{5} \begin{pmatrix} 1 & 6 & -9 & 6 \\ -8 & -13 & 12 & -3 \\ 13 & 8 & 3 & -12 \\ -6 & -1 & -6 & 9 \end{pmatrix} \sin(\delta)$$

$$1440 \begin{pmatrix} 4 & 3 & 0 & -3 \\ -5 & -2 & -3 & 6 \\ -2 & -5 & 6 & -3 \\ 3 & 4 & -3 & 0 \end{pmatrix} + 0\, L + 0 \sin(\delta)$$

```
H[1] = cCon[1] + cL[1] L + cs[1] s + cCross[1] L s;
Do[H[i] = cCon[i] + cL[i] L + cs[i] s, {i, {0, 2, 3}}]
```

`matrixK` $= K$; `raw` $= 3125\,(K - I)$ so $\frac{\text{raw}}{3125} + I = K$; `KCheck` is defined from the coefficient matrices $H$, and adjusted in the same way. All three agree,

```
eval = {σ → Sin[δ/2]^2, s → Sin[δ], L → 10., δ → (L + 1) θ};
KCheck = id + (4/3125 σ) Table[Table[H[k][[i, j]], {k, 0, 3}].σ^Range[0,3], {i, 4}, {j, 4}];
Chop[Max[Abs[raw/3125 + id - matrixK //. eval]]]
Chop[Max[Abs[KCheck - matrixK //. eval]]]
```

```
0

0
```

## The Gap Functions Derived and Analyzed

We seek an upper bound on the maximum absolute value of any entry in $K - I$. Start with $K - I$ and turn all coefficients positive.

```
raw0 = Expand[σ (H[0] + H[1] σ + H[2] σ^2 + H[3] σ^3)/3125] /. {r_Rational (σ^z_Integer s)/√5 :→ Abs[r] σ^z/√5,
    r_Rational σ^z_Integer :→ Abs[r] σ^z, r_Rational (σ s)/√5 :→ Abs[r] σ/√5, r_Rational σ :→ Abs[r] σ};
raw0[[1]]
```

$$\left\{\frac{6\,\sigma}{25}+\frac{3\,\sigma}{25\,\sqrt{5}}+\frac{6\,L\,\sigma}{25}+\frac{48\,\sigma^2}{25}+\frac{288\,\sigma^3}{625}+\frac{144\,\sigma^3}{625\,\sqrt{5}}+\frac{288\,L\,\sigma^3}{625}+\frac{1152\,\sigma^4}{625},\right.$$
$$\frac{21\,\sigma}{25}+\frac{18\,\sigma}{25\,\sqrt{5}}+\frac{6\,L\,\sigma}{25}+\frac{18\,\sigma^2}{25}+\frac{1008\,\sigma^3}{625}+\frac{864\,\sigma^3}{625\,\sqrt{5}}+\frac{288\,L\,\sigma^3}{625}+\frac{864\,\sigma^4}{625},$$
$$\left.\frac{6\,\sigma}{25}+\frac{27\,\sigma}{25\,\sqrt{5}}+\frac{6\,L\,\sigma}{25}+\frac{288\,\sigma^3}{625}+\frac{1296\,\sigma^3}{625\,\sqrt{5}}+\frac{288\,L\,\sigma^3}{625},\ \frac{9\,\sigma}{25}+\frac{18\,\sigma}{25\,\sqrt{5}}+\frac{6\,L\,\sigma}{25}+\frac{18\,\sigma^2}{25}+\frac{432\,\sigma^3}{625}+\frac{864\,\sigma^3}{625\,\sqrt{5}}+\frac{288\,L\,\sigma^3}{625}+\frac{864\,\sigma^4}{625}\right\}$$

Get everything in terms of $L$.

```
raw1 = Apart[Simplify[raw0 //. {σ → Sin[δ/2]^2, Sin[δ/2]^z_ :> δ^z/2^z, δ → 2 π/(√5 (L + 1))}]];
raw1[[1, 1]]
```

$$\left(\frac{6}{125\,(1+L)^8}+\frac{3}{125\,\sqrt{5}\,(1+L)^8}+\frac{42\,L}{125\,(1+L)^8}+\frac{18\,L}{125\,\sqrt{5}\,(1+L)^8}+\frac{126\,L^2}{125\,(1+L)^8}+\frac{9\,L^2}{25\,\sqrt{5}\,(1+L)^8}+\frac{42\,L^3}{25\,(1+L)^8}+\frac{12\,L^3}{25\,\sqrt{5}\,(1+L)^8}+\right.$$
$$\left.\frac{42\,L^4}{25\,(1+L)^8}+\frac{9\,L^4}{25\,\sqrt{5}\,(1+L)^8}+\frac{126\,L^5}{125\,(1+L)^8}+\frac{18\,L^5}{125\,\sqrt{5}\,(1+L)^8}+\frac{42\,L^6}{125\,(1+L)^8}+\frac{3\,L^6}{125\,\sqrt{5}\,(1+L)^8}+\frac{6\,L^7}{125\,(1+L)^8}\right)\pi^2+$$
$$\left(\frac{48}{625\,(1+L)^8}+\frac{192\,L}{625\,(1+L)^8}+\frac{288\,L^2}{625\,(1+L)^8}+\frac{192\,L^3}{625\,(1+L)^8}+\frac{48\,L^4}{625\,(1+L)^8}\right)\pi^4+$$
$$\left(\frac{288}{78125\,(1+L)^8}+\frac{144}{78125\,\sqrt{5}\,(1+L)^8}+\frac{864\,L}{78125\,(1+L)^8}+\frac{288\,L}{78125\,\sqrt{5}\,(1+L)^8}+\right.$$
$$\left.\frac{864\,L^2}{78125\,(1+L)^8}+\frac{144\,L^2}{78125\,\sqrt{5}\,(1+L)^8}+\frac{288\,L^3}{78125\,(1+L)^8}\right)\pi^6+\frac{1152\,\pi^8}{390625\,(1+L)^8}$$

All terms are like the preceding, with $L$ in the numerator and a power of $L+1$ in denominator. Next we turn $\frac{1}{L+1}$ into $\frac{1}{L}$, since we are chasing an upper bound.

```
raw2 = raw1 /. (1 + L)^q_ :> L^q;
raw2[[1, 1]]
```

$$\left(\frac{6}{125\,L^8}+\frac{3}{125\,\sqrt{5}\,L^8}+\frac{42}{125\,L^7}+\frac{18}{125\,\sqrt{5}\,L^7}+\frac{126}{125\,L^6}+\frac{9}{25\,\sqrt{5}\,L^6}+\frac{42}{25\,L^5}+\frac{12}{25\,\sqrt{5}\,L^5}+\frac{42}{25\,L^4}+\frac{9}{25\,\sqrt{5}\,L^4}+\right.$$
$$\left.\frac{126}{125\,L^3}+\frac{18}{125\,\sqrt{5}\,L^3}+\frac{42}{125\,L^2}+\frac{3}{125\,\sqrt{5}\,L^2}+\frac{6}{125\,L}\right)\pi^2+\left(\frac{48}{625\,L^8}+\frac{192}{625\,L^7}+\frac{288}{625\,L^6}+\frac{192}{625\,L^5}+\frac{48}{625\,L^4}\right)\pi^4+$$
$$\left(\frac{288}{78125\,L^8}+\frac{144}{78125\,\sqrt{5}\,L^8}+\frac{864}{78125\,L^7}+\frac{288}{78125\,\sqrt{5}\,L^7}+\frac{864}{78125\,L^6}+\frac{144}{78125\,\sqrt{5}\,L^6}+\frac{288}{78125\,L^5}\right)\pi^6+\frac{1152\,\pi^8}{390625\,L^8}$$

Then we use the assumption $L>\frac{2}{\epsilon}$ and we bound $\frac{1}{L}$ by $\frac{\epsilon}{2}$, which we accomplish by replacing $L$ by $\frac{2}{\epsilon}$.

```
raw3 = Expand[raw2 /. L → 2/ϵ];
raw3[[1, 1]]
```

$$\frac{3\,\pi^2\,\epsilon}{125}+\frac{21\,\pi^2\,\epsilon^2}{250}+\frac{3\,\pi^2\,\epsilon^2}{500\,\sqrt{5}}+\frac{63\,\pi^2\,\epsilon^3}{500}+\frac{9\,\pi^2\,\epsilon^3}{500\,\sqrt{5}}+\frac{21\,\pi^2\,\epsilon^4}{200}+\frac{9\,\pi^2\,\epsilon^4}{400\,\sqrt{5}}+\frac{3\,\pi^4\,\epsilon^4}{625}+\frac{21\,\pi^2\,\epsilon^5}{400}+$$
$$\frac{3\,\pi^2\,\epsilon^5}{200\,\sqrt{5}}+\frac{6\,\pi^4\,\epsilon^5}{625}+\frac{9\,\pi^6\,\epsilon^5}{78125}+\frac{63\,\pi^2\,\epsilon^6}{4000}+\frac{9\,\pi^2\,\epsilon^6}{1600\,\sqrt{5}}+\frac{9\,\pi^4\,\epsilon^6}{1250}+\frac{27\,\pi^6\,\epsilon^6}{156250}+\frac{9\,\pi^6\,\epsilon^6}{312500\,\sqrt{5}}+\frac{21\,\pi^2\,\epsilon^7}{8000}+\frac{9\,\pi^2\,\epsilon^7}{8000\,\sqrt{5}}+$$
$$\frac{3\,\pi^4\,\epsilon^7}{1250}+\frac{27\,\pi^6\,\epsilon^7}{312500}+\frac{9\,\pi^6\,\epsilon^7}{312500\,\sqrt{5}}+\frac{3\,\pi^2\,\epsilon^8}{16000}+\frac{3\,\pi^2\,\epsilon^8}{32000\,\sqrt{5}}+\frac{3\,\pi^4\,\epsilon^8}{10000}+\frac{9\,\pi^6\,\epsilon^8}{625000}+\frac{9\,\pi^6\,\epsilon^8}{1250000\,\sqrt{5}}+\frac{9\,\pi^8\,\epsilon^8}{781250}$$

Here is a quick approximation.

```
Collect[4/ϵ raw3[[1, 1]], ϵ] /. {zz_ ϵ^q_ :> N[zz] ϵ^q, zz_ ϵ :> N[zz] ϵ} /.
 {x_?NumericQ ϵ^q_ → Ceiling[100 x] /100. ϵ^q, x_?NumericQ π^q_ → Ceiling[100 x π^q] /100.}
```

$$0.95+3.42212\,\epsilon+5.3\,\epsilon^2+6.42\,\epsilon^3+6.53\,\epsilon^4+4.25\,\epsilon^5+1.45\,\epsilon^6+0.64\,\epsilon^7$$

The limit of the ratio to the desired $\frac{\epsilon}{4}$ is $\frac{12}{125}\pi^2$ which is less than 1. And if $\epsilon<\frac{1}{100}$, the $\epsilon$-expressions are all less than 0.9999.

But these expressions are increasing polynomials that vanish at 0, so they live in [0, 1] provided $0 < \epsilon < \frac{1}{100}$.

```
Max[Limit[raw3/(ϵ/4), ϵ → 0]]
N[Max[Limit[raw3/(ϵ/4), ϵ → 0]]]
```

$\frac{12\,\pi^2}{125}$

```
0.947482
```

```
Max[raw3/(ϵ/4) /. ϵ → 1/100.]
```

```
0.999898
```

And in fact one could replace the 2 by $\frac{24}{25}\,\frac{1}{5}\,\pi^2$, or about 1.895; the limit above will be 1 in this case.

## The Vanishing-Gap Subsequence for the Quadrahelix

The denominators of the convergents for $\frac{\theta}{2\,\pi}$, less 1, give the good $L$-values.

```
Denominator[Convergents[θ/(2. π)]] - 1
```

```
{0, 1, 2, 7, 10, 29, 40, 70, 182, 253, 1960, 12 019, 13 980, 26 000,
 143 985, 601 944, 5 561 490, 6 163 435, 11 724 926, 17 888 362, 65 390 015}
```

## The Limiting Angles

The dot product gives the cosine, the four angles sum to $2\,\pi$, and alternate angles are equal, so this tells the whole story.

```
L1 = TrigExpand[ExpandAll[base2[[1]] - base1[[2]]]];
L2 = TrigExpand[ExpandAll[TrigExpand[base3[[1]] - base2[[1]]]]];
dot = L1.L2;
Limit[dot/(Sqrt[L1.L1] Sqrt[L2.L2]) /. L θ → -θ, L → ∞]
```

$-\frac{1}{5}$

## The Limiting Rhombus

To get the rhombus we can take the first two points to be (0, 0) and (0, 1).

```
Clear[A]; A[n_] := {x[n], y[n]};
x[1] = y[1] = 0; x[2] = 0; y[2] = 1;
eqnsUnitLength = Table[1 == (A[Mod[i + 1, 4, 1]] - A[i]).(A[Mod[i + 1, 4, 1]] - A[i]), {i, 4}];
eqnsAngle = Table[(-1)^i/5 == (A[Mod[i + 2, 4, 1]] - A[Mod[i + 1, 4, 1]]).(A[Mod[i + 1, 4, 1]] - A[Mod[i, 4, 1]]), {i, 4}];
A /@ Range[4] /. Solve[Join[eqnsUnitLength, eqnsAngle], Join[A[3], A[4]]]
```

$\left\{\left\{\{0, 0\}, \{0, 1\}, \left\{-\frac{2\sqrt{6}}{5}, \frac{4}{5}\right\}, \left\{-\frac{2\sqrt{6}}{5}, -\frac{1}{5}\right\}\right\}, \left\{\{0, 0\}, \{0, 1\}, \left\{\frac{2\sqrt{6}}{5}, \frac{4}{5}\right\}, \left\{\frac{2\sqrt{6}}{5}, -\frac{1}{5}\right\}\right\}\right\}$

# Diophantine Approximation by LatticeReduction

```
DiophantineApproximationBabaiLLL[{α_, β_, γ_}, X_] := Module[{s, a, b, c, B, BStar, t, x, y},
  ϵ = (3 Abs@β / X)^2; s = 1/ϵ^2 Max[Abs[{1, α, β}]];
  {a, b, c} = Round[s {α, β, -γ}];
  B = LatticeReduce[{{a, 1, 0}, {b, 0, 1}}];
  BStar = {B⟦1⟧, B⟦2⟧ - (B⟦2⟧.B⟦1⟧ B⟦1⟧)/(B⟦1⟧.B⟦1⟧)};
  t = {c, 0, 0};
  Do[t -= Round[t.BStar⟦j⟧ / (BStar⟦j⟧.BStar⟦j⟧)] B⟦j⟧, {j, 2, 1, -1}];
  {x, y} = t⟦{2, 3}⟧;
  {x, y, N[Log10[Abs[x α + y β - γ]], 20]}];

Grid[Prepend[Table[DiophantineApproximationBabaiLLLA[{θ, 2 π, γ}, 10^i], {i, 2, 6.5, 1/2}] /.
   {x_Integer, y_, z_} :> {x, y, NumberForm[z, 3]} /. {x_Integer /; x < 0, y_, z_} :>
   {-x - 1, -y + 1, NumberForm[eL = N[Log10[Abs[(-x - 1) θ + (-y + 1) 2 π - γ1]], 20], 3]}, {"x", "y", "log10error"}]]
```

| x | y | $\log_{10}$error |
|---|---|---|
| 4 | -1 | -1.80 |
| 686 | -251 | -3.46 |
| 1274 | -466 | -3.88 |
| 64 708 | -23 692 | -5.48 |
| 666 653 | -244 088 | -5.65 |
| 1 870 543 | -684 880 | -6.86 |
| 111 021 125 | -40 649 248 | -7.79 |
| 233 817 317 | -85 609 817 | -9.23 |
| 3 113 400 370 | -1 139 939 675 | -9.71 |
| 434 337 601 428 | -159 028 266 709 | -10.4 |